\input amstex
\documentstyle{amsppt}
\refstyle {A}
\magnification=\magstep1 

\hoffset .25 true in
\voffset .2 true in

\hsize=6.1 true in
\vsize=8.5 true in

\def\a{\alpha}
\def\b{\beta}

\def\d{\delta}
\def\l{\lambda}
\def\m{\mu}
\def\n{\nu}

\def\om{\omega}

\def\G{\Gamma}
\def\D{\Delta}
\def\HD{\hat\Delta}
\def\L{\Lambda}
\def\Th{\Theta}

\def\A{\bold A}
\def\B{\bold B}
\def\CB{\Cal B}

\def\F{\Bbb F}
\def\FF{\bold F}

\def\M{\Cal M_c}
\def\PP{\Cal P}
\def\Q{\operatorname{Q}}
\def\R{\Bbb R}
\def\RR{\Cal R}
\def\SS{\Cal S}
\def\S{\bold S}
\def\T{\bold T}
\def\TT{\Cal T}

\def\op{\operatorname}
\def\lrar{\leftrightarrow}
\def\ov{\overline}
\def\sub{\subset}
\def\subeq{\subseteq}
\def\Id{\operatorname{Id}}

\def\bu{\bullet}
\def\hx{\hat X}
\def\hk{\hat K}
\def\hook{\hookrightarrow}
\def\Hom{\operatorname{Hom}}
\def\Vect{\operatorname{Vect}_{\Bbb F}}
\def\supp{\operatorname{supp}}
\def\refn{\operatorname{ref}}

\def\vplim{\varprojlim}

\def\ti{\tilde}
\def\wti{\widetilde}
\def\php{\Phi^+}
\def\phm{\Phi^-}

\def\psp{\Psi^+}
\def\psm{\Psi^-}
\def\upp{\Upsilon^+}
\def\upm{\Upsilon^-}

\def\SH{\Cal {SH}}
\def\CC{{\Cal C}^b}
\def\DD{{\Cal D}^b}
\def\moda{\operatorname{mod -}\! A}
\def\modb{\operatorname{mod -}\! B}
\def\modc{\operatorname{mod -}\! C}
\def\modgra{\operatorname{modgr -}\! A}
\def\modgrb{\operatorname{modgr -}\! B}
\def\modgrc{\operatorname{modgr -}\! C}

\def\HOM{\bold{Hom}}
\def\I{\bold I}
\def\codim{\operatorname{codim}}
\def\Pre{\operatorname{Pre}}
\def\LL{\bold L}
\def\II{\Cal I}
\def\U{\bold U}
\def\V{\bold V}

\def\PD{{}^{-\delta}\!\Delta}

\def\sup{\supset}
\def\supeq{\supseteq}

\def\ver{\text{\sl Vert}}
\def\e{\varepsilon}
\def\vp{\varphi}
\def\ves{\emptyset}


\topmatter
\title
Sheaves on Triangulated Spaces and Koszul Duality
\endtitle
\author
Maxim Vybornov
\endauthor
\address Duke University,
Department of Mathematics,
Durham, NC 27708
\endaddress
\email mv\@math.duke.edu
\endemail
\dedicatory
Dedicated to Robert MacPherson on the occasion of his 55th birthday
\enddedicatory
\abstract
Let $X$ be a finite connected simplicial complex,
and let $\delta$ be a perversity 
(i.e., some function from integers to integers).
One can consider two categories: (1) the
category of perverse sheaves cohomologically constructible with
respect to the triangulation, and (2) the category of sheaves constant
along the perverse simplices ($\delta$-sheaves). We interpret
the categories (1) and (2) as categories of modules over 
certain quadratic (and even Koszul) algebras $A(X,\delta)$
and $B(X,\delta)$ respectively, and we prove that $A(X,\delta)$
and $B(X,\delta)$ are Koszul dual to each other.
We define the $\delta$-perverse 
topology on $X$ and prove that the category of sheaves
on perverse topology is equivalent to the category of  
$\delta$ sheaves.
Finally, we study the relationship between the Koszul
duality functor and the 
Verdier duality functor
for simplicial sheaves and cosheaves.
\endabstract

\endtopmatter


\subhead
{Introduction}
\endsubhead

\subhead{1}\endsubhead
The study of constructible sheaves
on a cell complex $X$ leads to the notion of a {\sl cellular sheaf},
which was developed by W. Fulton, M. Goresky, R. MacPherson,
and C. McCrory in a seminar at Brown University in 1977-78.
The systematic exposition of the theory of
cellular sheaves has been presented in the
A. Shepard's Doctoral Thesis [Shep] (cf. [Kash]). 

A cellular sheaf is a gadget which assigns vector spaces
to cells in $X$ and linear maps to pairs of incident cells.
(Cellular sheaves can also be interpreted as sheaves on the 
finite topology generated by open stars of cells.)
It is easy to interpret such linear algebra gadgets as modules
over an associative algebra $B(X)$.

In this paper we will work with a 
finite connected simplicial complex $X$. 
We can consider simplicial complexes without the loss of generality
since any reasonable stratified space can be triangulated [Gor].
The category
of constructible sheaves of $\F$-vector spaces on $X$ is denoted by 
$\SH_c(X)$.
We formulate here the basic result of the cellular sheaf theory.

\proclaim{Theorem A} The following categories
$$
\SH_c(X)\simeq\modb(X)   
$$
are equivalent.
\endproclaim

\subhead{2}\endsubhead
A {\sl perversity function} $\ov p:\Bbb Z\to\Bbb Z$ is
a key ingredient in the definition of the intersection homology
groups introduced by M. Goresky and R. MacPherson in [GM1].
This seminal paper also contains the construction of the 
{\sl basic sets}. The later version of basic sets 
called {\sl perverse skeleta}, and their building blocks
called perverse cells were introduced in [Mac2],
as well as the cellular version $\d:\Bbb Z\to\Bbb Z$ of
the perversity function. (The classical and cellular perversity 
functions are in one-to-one correspondence.)
An extremely important application of the idea of perverse cell 
complex appears in the work of I. Mirkovic and K. Vilonen [MkVi]. 
They use perverse cells in loop Grassmannians associated
to complex algebraic groups to give a proof of the geometric
Satake isomorphism theorem announced by V. Ginzburg in [Ginz].

Given a cellular perversity function $\d$ we construct the
$\d$-perverse finite topology on $X$.
The idea of the existence 
of perverse topology appeared in the discussions 
surrounding the lectures [Mac2], and it was communicated to us
by M. Goresky and R. MacPherson. 
One recovers the topology generated by open stars
for {\sl bottom} (i.e. $\d(k)=-k$, $k\geq 0$) cellular perversity. 
Sheaves on perverse topology can be interpreted as modules
over an associative algebra $B(X,\d)$. In other words, we have

\proclaim{Theorem B1}The following categories
$$
\SH(\TT(X,\d))\simeq\modb(X,\d)   
$$
are equivalent.  
\endproclaim 

Let us consider the category $\SH_c(X,\d)$ of sheaves on $X$
constant along $(-\d)$-perverse simplices ($\d$-sheaves).
We give a linear algebra description of the category of 
$\d$-sheaves, generalizing Theorem A.
Our linear algebra
gadgets assign vector spaces to perverse simplices and linear maps to
pairs of ``incident'' perverse simplices.
As in the classical case, 
it is easy to interpret such linear algebra gadgets as modules
over the algebra $B(X,\d)$. In other words, we have

\proclaim{Theorem B2} The following categories
$$
\SH_c(X,\d)\simeq\modb(X,\d)   
$$
are equivalent.
\endproclaim

We believe that it is very important to study $\d$-sheaves
for more general stratified spaces. 
The construction of $\d$-sheaves on flag varieties
should yield new geometric
realization of various representation theoretic constructions.
In particular, for a parabolic flag variety this approach could 
lead to a geometric realization of a singular block of the 
Bernstein-Gelfand-Gelfand category $\Cal O$,
which would be useful for a geometric
categorification of the Temperley-Lieb algebra (see [BFK]).

\subhead{3}\endsubhead
Let $X$ be a finite connected simplicial complex,
and let $\d$ be a perversity.
The category of cohomologically constructible with respect to the 
triangulation perverse sheaves on $X$ [BBD, Mac 1] is denoted by $\M(X,\d)$.
A lot of work has been done to represent perverse sheaves
in terms of linear algebra data (e.g. [MaVi], [BrGr]).
The linear algebra data
description of the category $\M(X,\d)$  was obtained by 
R. MacPherson in [Mac2, Mac3]. 
MacPherson's  linear algebra gadgets are called 
{\sl cellular perverse sheaves}.
It is easy to interpret cellular perverse sheaves
as modules over an associative algebra $A(X,\d)$.
In other words, we have

\proclaim{Theorem C} The following categories
$$
\M(X,\d)\simeq\moda(X,\d)
$$
are equivalent.
\endproclaim

\subhead{4}\endsubhead
In the present and the following sections we describe our main 
new results: Koszul duality patterns for sheaves
and applications. The proof of the Theorem D2 was sketched in [Vyb2].
Theorem D1, Corollary D3, and Theorem E are new.

The notion of a {\sl Koszul} algebra was introduced by
S. Priddy [Prid]. In many important cases the algebras 
underlying the categories of
perverse sheaves turn out to be Koszul (cf. [BGSc, BGSo, PS]). 
The Koszulity of $A(X,\d)$ and $B(X,\d)$ has been 
proved by A. Polishchuk in [Pol].
Given a Koszul category, it seems
to be very important to investigate the Koszul dual category
and the duality functor.
We establish the following claims:

\proclaim{Theorem D1} Perverse sheaves on the
ordinary topology of $X$ are {\rm Koszul dual} to the ordinary
sheaves on the perverse topology. 
(For the precise formulation see 4.2.11.)
\endproclaim

\proclaim{Theorem D2} Perverse sheaves constant
along ordinary simplices of $X$ are {\rm Koszul dual} to ordinary
sheaves constant along perverse simplices.
(For the precise formulation see 4.2.11.)
\endproclaim  

The basic notion underlying the definition of the Intersection Homology 
is that
of the {\sl allowable} set (see [GM1] and 4.2.12).

\proclaim{Corollary D3} The Koszul duality functor transforms
a stratification-constructible perverse sheaf on a pseudomanifold into 
a complex of sheaves with allowable support.    
(For the precise formulation see 4.2.12.)
\endproclaim

One of the properties of the Koszul duality functor is
that it preserves the hypercohomology of a complex of 
sheaves.

If one considers a perverse sheaf as 
a ``generalized homology theory'' (F\'ary functor
[Mac4, GMMV]), then Theorem D2 gives a way
to see ``geometric chains'' of this homology theory.

\subhead{5}\endsubhead
Let $D$ be the Verdier duality functor acting
in the derived category of simplicial sheaves.
We introduce another functor $\wti{D}$, which we call
the ``Verdier duality for cosheaves.'' We prove that
the Koszul duality functor $L$ intertwines $D$ and $\wti{D}$.  
More precisely, we have

\proclaim{Theorem E} The following functors
$$
L\circ D\simeq\wti{D}\circ L
$$
are isomorphic.
\endproclaim

The question of how Koszul duality commutes with Verdier
duality motivating Theorem E was formulated by V. Ginzburg
in connection with Beilinson-Ginzburg conjecture 
[BG, 5.18, 5.24].

\subhead{6}\endsubhead
We try to give the explicit reference in the text each time
we borrow some material. However, we would like to
emphasize that our main sources are 
[BBD], [BGSo], [CPS], [GM1], [KS], [Mac2, Mac3], [Pol], and [Shep].
Some standard fact listed in the text
without reference is probably lifted from one of the above.

\subhead{7}\endsubhead
The paper is organized as follows.
In Chapter 1 we list some standard facts mostly lifted
from classical books and papers [BBD, GeMa, GM1, GM2, KS, Mas] and 
other sources, in particular [Mac2], and more recent papers
[Pol, Vyb1, Vyb2].
In Chapter 2 we study the relationship between sheaves and
presheaves on a base of topology. We prove Theorem A, Theorem B1,
and list some properties of simplicial sheaves.
Chapter 3 is devoted to sheaves constant along perverse simplices.
We prove Theorem B2 and study some properties of $\d$-sheaves.
In Chapter 4 we are initially concerned with the perverse algebras. 
Then we study the Koszul duality functors for sheaves
and their properties. We complete the proof of Theorem C,
Theorem D1, Theorem D2, Corollary D3, and Theorem E.

\subhead{8}\endsubhead
We list here some notational conventions adopted in the paper. 
Unless specified otherwise, $X$ is assumed
to be a finite connected simplicial complex.
Other conventions: 
$\Bbb Z$ integer numbers;
$\R$ real numbers;
$\F$ commutative field of $\op{char}=0$;
$\FF$ constant sheaf associated to $\F$;
$\Id$ identity operator (morphism, functor).
All vector spaces are assumed to be over $\F$ unless
specified otherwise.

\subhead{Acknowledgement}\endsubhead
This paper is based on my Doctoral Dissertation [Vyb3].

I am deeply grateful to my advisor Igor Frenkel for
countless hours he spent teaching me beautiful and exciting 
mathematics, and his relentless encouragement and support
during my four years at Yale.

My work is based on the theory of cellular perverse
sheaves discovered by Robert MacPherson. I had a
privilege to attend his wonderful M.I.T. lectures on the subject. 
Without these lectures, discussions with, and all kinds of support 
from
Robert MacPherson the present work would have never been 
started, let alone finished.

I am greatly indebted to Mark Goresky who supported and helped 
to shape this work
on every single stage from the very beginning to the very end. 
I am sincerely appreciating his continuous
mathematical and professional advice and encouragement, 
and his sense of humor, for many years.
I have greatly benefited from numerous enlightening discussions with
A. Beilinson, V. Ginzburg, M. McConnell, I. Mirkovic,
A. Polishchuk, L. Positselski, C. Procesi, K. Vilonen,
and G. Zuckerman.
I am grateful to R. Lee for inviting
me to speak on the topic at Yale Topology Seminar, and to
M. Khovanov for inviting me to speak on the topic
at minicourse on perverse sheaves organized by him at the
Institute for Advanced Study.
I would also like to express my sincere appreciation of the
warm and stimulating atmosphere created by the faculty
and staff, and my fellow graduate students  
at the Mathematics Department at Yale. I am especially
indebted to A. Malkin, Y.-H. Kiem, A. Liakhovskaia, I. Pak, and 
M. Yampolsky. 
I am also grateful to the Institute for Advanced Study
for their hospitality in June 1995, and in the Fall of 1998.
I am grateful to Yale University for the generous  
financial support and to Alfred P. Sloan Foundation
which selected me as Doctoral Dissertation Fellow
in 1998-99. I have also been financially 
supported by the National Science Foundation on various occasions
starting from the summer of 1995.
Finally, I would like to thank the anonymous referee for useful remarks.

\head{Chapter 1. Preliminaries}\endhead

\subhead 1.1. Simplicial complexes and perversities\endsubhead
The bulk of this section is standard algebraic topology
lifted from [GeMa, KS, Mas]. The material on perversities
is borrowed from [BBD, GM1, Mac2, Pol, Vyb2].

\subhead{1.1.1}\endsubhead 
Let $\D=\{v_0,v_1,\dots v_r\}$ be the $r$-dimensional 
nondegenerate simplex spanned by vertices $\{v_0,v_1,\dots v_r\}$.
The barycenter of $\D$ will usually be denoted by $c$.

\subhead{1.1.2}\endsubhead 
Let $K$ be a simplicial set with the set of vertices $\ver$. 
Let $\R^{\ver}$ denote the set of maps from $\ver$ to $\R$. An element
$x\in\R^{\ver}$ is nothing but a family $x(v)\in\R$, indexed by 
$v\in\ver$. We equip $\R^{\ver}$ with the usual Euclidean topology.
To a simplex of $K$ we associate the geometric simplex $\D$ in $\R^{\ver}$ by:
$$
\D=\{x\in\R^{\ver} \ :\  x(v)=0 \text{ for } v\not\in\D,
x(v)>0 \text{ for } v\in\D \text{ and } \sum_{v}x(v)=1\}.
$$ 
The geometric realization of $K$ denoted by $|K|$ is called a 
flat (or Euclidean) simplicial complex. 
A finite {\sl simplicial complex} $X$ is a topological space
homeomorphic to $|K|$ for some simplicial set $K$.
The topology induced from $|K|$ will sometimes be refered to
as usual topology. 
A simplex $\D$ of $X$ is a homeomorphic image of the simplex $\D$ 
of $|K|$. We say that two simplices $\D$ and $\D'$ are incident
if either $\D$ is a face of $\D'$ or $\D'$ is a face of $\D$.
If $\D$ and $\D'$ are incident, we write $\D\lrar\D'$.
Unless specified otherwise,  we will always assume that $X$ 
is a finite connected simplicial complex.
We assume that $\dim X=n$ (i.e. X has some $n$-dimensional simplices
and does not have any $(n+1)$-dimensional simplices). 
 
\subhead {1.1.3}\endsubhead
Since $X$ is a (regular) cell complex we can choose an
orientation for each simplex (see [Mas, IX.5]) and define the 
{\sl incidence numbers} $[\D:\D']$ for any pair of simplices such 
that $\D'$ is a $\op{codim}\ 1$ face of $\D$. We record
the following

\proclaim{Lemma [Mas, IX.7.1]} Let $\D',\D$ be as above. Then  
$[\D:\D']=\pm 1$.
\endproclaim

\subhead{1.1.4}\endsubhead
Let $K$ be a finite simplicial set.
The first barycentric subdivision of $K$ is denoted by $\hk$.
If $X$ is a simplicial complex then its  
first barycentric subdivision is denoted by $\hx$.

\definition{1.1.5. Definition} 
\roster
\item ([Mac2]) A {\sl cellular perversity $\d :{\Bbb Z}_{\geq 0}\to 
{\Bbb Z}$} is a function from the
non-negative integers ${\Bbb Z}_{\geq 0}$ to the integers such that 
$\delta (0)=0$ and $\delta$ takes every interval  
$\{0,1,\ldots ,k\}\subset {\Bbb Z}_{\geq 0}$ bijectively to an interval
$\{a,a+1,\ldots ,a+k\}\subset {\Bbb Z}$ for some 
$a\in {\Bbb Z}_{\leq 0}$.
In other words, a perversity is such a function $\delta$ that 
$\delta(0)=0$ and for $k\in\Bbb Z_{\geq 0}$,
$$\delta(k)=
\cases \text {either} & \underset{i\in[0,k-1]}\to{\max} \delta(i)+1 \\ 
\text { or} & \underset{i\in[0,k-1]}\to{\min} \delta(i)-1. \\ \endcases 
$$
\item (cf. [GM1]) A {\sl classical} (or {\sl Goresky-MacPherson})
{\sl  perversity} $\ov p :{\Bbb Z}_{\geq 0}\to 
{\Bbb Z}_{\geq 0}$ is a monotonously non-decreasing function 
from the set of
non-negative integers ${\Bbb Z}_{\geq 0}$ to itself such that 
$\ov p (0)=0$ and 
$\ov p(k)-\ov p(k-1)$ is either $0$ or $1$. 
\item ([BBD]) A {\sl BBDG perversity} $p :{\Bbb Z}_{\geq 0}\to 
{\Bbb Z}_{\leq 0}$ is a function such that $p(0)=0$,
$0\leq p(m)-p(n)\leq n-m$ for $m\leq n$. 
We define the {\sl dual perversity} by $p^*(k)=-p(k)-k$.
\endroster  
\enddefinition

\example{Example} We will often work with two extreme cellular 
perversities, {\sl top perversity} $\d(k)=k$, $k\geq 0$,
and {\sl bottom perversity} $\d(k)=-k$, $k\geq 0$.
\endexample

Let $p$ be a BBDG perversity. An integer $k\in\Bbb Z_{\geq 0}$ 
is of {\sl type} $*$ (resp. {\sl type} $!$) if
$p(k)=p(k-1)-1$ (resp. $p(k)=p(k-1)$).
It is convenient to assume that $0$ is of
both types (see [Pol, 1]).

\subhead 1.1.6 \endsubhead In this subsection we study the
relationship between different versions of perversity. 
It is easy to see that if $\ov p$ is a classical perversity, then 
$p=-\ov p$ is a BBDG perversity and vice versa. 

\proclaim{Lemma} There is a one-to-one correspondence 
between cellular perversities and BBDG perversities.
\endproclaim

\demo{Proof} Let $p$ be a BBDG perversity.
The corresponding cellular perversity $\d$ is defined by:
$$
\d(k)=\cases
-p(k), & k \text{ of type } *\ \\
p^*(k), & k \text{ of type } !\ .\\
\endcases
$$
\enddemo

Below we work mostly with the cellular version 
of perversity, but sometimes we have to refer to the associated
classical and BBDG perversity since their use is standard
in the context of Intersection Homology and perverse sheaves
respectively.

\subhead 1.1.7 \endsubhead 
Until the end of this section 
$X$ is assumed to be a finite connected simplicial complex, 
$\dim X=n$.
In this subsection we recall the definition of the basic sets
given in [GM1]. 
Let us fix a classical perversity $\ov p$.
For an integer $i\geq 0$ define the function $L_i^{\ov p}$ as 
follows:
$$L_i^{\ov p} (0)=i, \qquad L_i^{\ov p} (n+1)=-1,$$
and if $1\leq c\leq n$ set:
$$
L_i^{\ov p} (c)=\cases -1 & \text{ if } i-c+\ov p(c)\leq -1 \\
n-c & \text{ if } i-c+\ov p(c)\geq n-c \\ 
i-c+\ov p(c) & \text{ otherwise. } \\ 
\endcases
$$
Define $dL_i^{\ov p} (c)=L_i^{\ov p} (c)-L_i^{\ov p} (c+1)$ 
(which is either $0$ or $1$).
Define $Q_i^{\ov p}$ to be the subcomplex of 
$\hx$ spanned by the set of barycenters
of such simplices $\D$ of $X$ that: 
$$d L_i^{\ov p} (n-\dim\D)=1.$$

\definition{1.1.8. Definition [Mac2]}
Given a perversity $\d$, 
we define $\d(\D)=\d(\dim\D)$, where $\D$ is a simplex of $X$.
Given a simplex $\HD=\{c_0,c_1,\dots,c_s\}$ of $\hx$ (where 
$\{c_0,c_1,\dots,c_s\}$ are barycenters of 
$\D_0, \D_1,\dots,\D_s$) 
we denote by $\max\HD$ such
vertex $c_i$ that $\d(\D_i)=\max\{\d(\D_0),\d(\D_1),\dots, \d(\D_s)\}$.
Given a simplex $\D$ with the barycenter $c$
we define the corresponding {\sl perverse simplex}:
$$
^\d\!\D=\bigsqcup_{\max\HD=c} \HD.
$$
We define the {\sl $k$-th perverse skeleton} $X^{\d}_k$, 
$\min_{[0,n]}\d\leq k\leq\max_{[0,n]}\d$, as follows:
$$
X^{\d}_k=\bigsqcup_{\d(\D)\leq k}{}^\d\!\D\sub X.
$$
Thus, we have a filtration:
$$
X^{\d}_i\sub X^{\d}_{i+1}\sub\dots\sub X^{\d}_{i+n-1}
\sub X^{\d}_{i+n}=X,
$$
where $i=\min_{[0,n]}\d$.
It is easy to see that perverse simplices are connected 
components of $X^{\d}_k-X^{\d}_{k-1}$, 
$\min_{[0,n]}\d\leq k\leq\max_{[0,n]}\d$. The decomposition of 
$X$ into a disjoint union of perverse simplices is called
{\sl $\d$-perverse triangulation} of $X$.
\enddefinition

\proclaim{1.1.9. Lemma} Let $\d$ be the cellular perversity 
corresponding to a classical perversity $\ov p$. 
Let $\ov q(k)=\ov p(n)-\ov p(n-k)$ be another perversity. Then
$$
Q_i^{\ov q}=X^{\d}_{\ov p(n)-n+i}\ .
$$
\endproclaim

\demo{Proof} The cellular perversity $\d$
corresponding to a classical perversity $\ov p$ is given by:
$$
\d(k)=\cases
\ov p(k), & k \text{ of type } *\ \\
\ov p(k)-k, & k \text{ of type } !\ .\\
\endcases
$$
Here $k$ is of type $*$ (resp. type $!$) if 
$\ov p(k)=\ov p(k-1)+1$ (resp. $\ov p(k)=\ov p(k-1)$).
Let $m$ be an integer $m\in\Bbb Z_{\geq 0}$.
It is easy to see that the following conditions are
equivalent:
\roster
\item $d L_i^{\ov q} (n-m)=1$.
\item $\d(m)\leq\ov p(n)-n+i$.
\endroster 
\enddemo

\subhead 1.2. Partially ordered sets, presheaves, and  
algebras \endsubhead
In this section we introduce the notion of a presheaf 
on a partially ordered set, and the notion of a
simplicial perverse sheaf. We notice that presheaves
and simplicial perverse sheaves can be interpreted as modules
over certain quadratic algebras.

\subhead {1.2.1}\endsubhead
Let $\L$ be a partially ordered set with the relation $\a\leq\b$.
We can associate to $\L$ a small category whose objects are
elements of $\L$ and the set of morphisms $\Hom(\a,\b)$ contains
one element if $\a\leq\b$, and is empty otherwise. By abuse of
notation we will also denote this category by $\L$. 

\definition{Definition} A presheaf $\S$ on $\L$ consists of a vector
space $S(\a)$ associated to every $\a\in\L$ (stalk of $\S$ at $\a$),
and a linear map $s(\a,\b):S(\a)\to S(\b)$, for $\a\geq\b$
(restriction map) asssociated to every $\a\geq\b$ in such a way that
$\S$ is a functor from the category opposite to $\L$ to the category of
vector spaces. A morphism $f$ between two presheaves $\S$ and $\S'$ is 
a collection of stalkwise linear maps commuting with the restriction maps.
The category of presheaves on $\L$ is denoted by $\Pre\L$.
\enddefinition
\subhead {1.2.2}\endsubhead
Until the end of this section 
$X$ is assumed to be a finite connected simplicial complex.
Let us fix a perversity $\d$. We will introduce a partial
order on the set of simplices of $X$ as follows (cf. [Pol, 2]). 
We say that 
$\D\geq\D'$ if there exists a sequence of simplices
$\D=\D_0, \D_1, \D_2,\dots,\D_r=\D'$ such that:
\roster
\item $\D_i\lrar\D_{i+1}\lrar\D_{i+2}$, $0\leq i\leq r-2$,
\item $\d(\D_i)=\d(\D_{i+1})+1$, $0\leq i\leq r-1$.
\endroster
The set of all simplices of $X$ with this partial order is 
denoted by $\L(X,\d)$, and from now on $\D\geq\D'$ will
refer to the relation in $\L(X,\d)$.

\proclaim{1.2.3. Lemma} (a) Let $\D\lrar\D'$, $\d(\D)\geq\d(\D')$.
Then $\D\geq\D'$. Moreover, we can choose a sequence
of simplices such that $\D\lrar\D_i\lrar\D'$ for all $i$.

(b) If $\D\geq\D'$ and $\d(\D)\geq\d(\D')\geq 0$ or
$0\geq\d(\D)\geq\d(\D')$, then $\D\lrar\D'$.

(c) If $\D>\D'$ and $\D\not\lrar\D'$, then $\d(\D)>0$, 
$\d(\D')<0$. If in addition $\d(\D)-\d(\D')=2$, then
$\d(\D)=1$, $\d(\D')=-1$.  
\endproclaim

\demo{Proof} The proof is elementary and is left to the reader.
\enddemo

\subhead {1.2.4}\endsubhead
We define a category $\RR (X,\d)$ of linear algebra data.
\definition{Definition} 
An object $\bold S$ of $\RR (X,\d)$
is the following data:
\roster
\item  (Stalks) For every simplex $\D$ in $X$, a finite dimensional
vector space 
$S(\D)$ called the {\sl stalk} of ${\S}$ at $\D$,
\item  (Restriction  maps) For every pair of simplices $\D$ and $\D'$ in 
$X$ such that $\d(\D)=\d(\D')+1$, and $\D\lrar\D'$, a linear map 
$s(\D ,\D'):S(\D)\to S(\D')$ called the 
{\sl restriction map},
\endroster
subject to the following ``equivalence'' axiom:  

Suppose that $\D'$, $\D''$, $\D_1$ and $\D_2$ are simplices
of $X$  such that:
\roster
\item $\d(\D')=k+1$, $\d (\D'')=k-1$, and $\d(\D_1)=\d(\D_2)=k$, 
\item $\D'\lrar\D_1\lrar\D''$ and $\D'\lrar\D_2\lrar\D''$.
\endroster
Then
$$
s(\D_1, \D'')\circ s(\D',\D_1)=
s(\D_2, \D'')\circ s(\D',\D_2).
$$
A morphism $h$ in this category is the collection of 
stalkwise linear maps $h(\D)$ commuting with
the restriction maps (cf. 1.2.1).
\enddefinition

\subhead{1.2.5}\endsubhead
In this subsection $n$ is not supposed to be the dimension of $X$. 

\proclaim{Lemma-Definition} Let $\D_I\geq\D_T$ be two simplices, 
and let
$\S\in\RR(X,\d)$. We define the linear map
$s(\D_I,\D_T):S(\D_I)\to S(\D_T)$ as follows.
If $\D_I=\D_T$, then $s(\D_I,\D_T)=\Id_{S(\D_I)}$.
Let $\D_I>\D_T$, and let  
$$
\D_I=\D_0 >\D_1>\D_2>\dots > \D_n=\D_T\ ,\tag {$*$}
$$
$\d(\D_i)=\d(\D_{i+1})+1$, $i=0, \dots, n-1$; and
$$
\D_I=\D'_0 >\D'_1>\D'_2>\dots > \D'_n=\D_T\ ,\tag {$**$}
$$
$\d(\D'_i)=\d(\D'_{i+1})+1$, $i=0, \dots, n-1$ 
be two sequences of simplices of $X$. Then
$$
s(\D_I,\D_T):=s(\D_{n-1},\D_n)\circ\dots\circ s(\D_0,\D_1)
=s(\D'_{n-1},\D'_n)\circ\dots\circ s(\D'_0,\D'_1).\tag\dag$$
\endproclaim

\demo{Proof} The proof is elementary and is left to the reader. 
\enddemo

\proclaim{1.2.6. Lemma} The following categories
$$
\Pre\L(X,\d)=\RR(X,\d)
$$
are isomorphic.
\endproclaim

\demo{Proof} We leave it to the reader to define  
the functor $F:\Pre\L(X,\d)\to\RR(X,\d)$, 
the functor $G:\RR(X,\d)\to\Pre\L(X,\d)$,
and verify that $F\circ G=\Id$,
and $G\circ F=\Id$.
\enddemo

If $\d=\text{bottom perversity}$ 
(resp. $\d=\text{top perversity}$),
then objects of $\Pre\L(X,\d)=\RR(X,\d)$
are called {\sl simplicial sheaves}
(resp. {\sl simplicial cosheaves}).

\subhead {1.2.7} \endsubhead
We will introduce another category $\PP(X,\d)$ 
of linear algebra data.

\definition{Definition [Mac2]}
An object $\bold S$ of $\PP (X,\d)$
is the following data:
\roster
\item  (Stalks) For every simplex $\D$ in $X$, a finite dimensional
vector space 
$S(\D)$ called the {\sl stalk} of ${\S}$ at $\D$,
\item  (Boundary  maps) For every pair of simplices $\D$ and $\D'$ in 
$X$ such that $\d(\D)=\d(\D')+1$, and $\D\lrar\D'$, a linear map 
$s(\D ,\D'):S(\D)\to S(\D')$ called the 
{\sl boundary map},
\endroster
subject to the following ``chain complex'' axiom:  

If $\D'$ is any simplex such that
$\d(\D')=k+1$ and $\D''$ is any simplex such that 
$\d(\D'')=k-1$, then
$$\sum\limits_{\matrix
\D : \d(\D)=k,\\
\D'\lrar\D\lrar\D''\endmatrix}
s(\D ,\D'')\circ s(\D',\D )=0.$$
The morphisms in this category are 
stalkwise linear maps commuting with
the boundary maps. Objects of $\PP(X,\d)$ 
are called {\sl simplicial perverse sheaves}.
\enddefinition

\subhead {1.2.8}\endsubhead
A quiver $\Q$ is a finite simple oriented tree.
We will denote the set of vertices of $\Q$ by $V(\Q)$
and the set of arrows by $E(\Q)$.

We associate a quiver $\Q(X,\d)$ to the set $\L(X,\d)$ as follows.
The vertices of $\Q(X,\d)$ are indexed by the elements
of $\L(X,\d)$ (i.e. $V(\Q(X,\d))=\L(X,\d)$), and there is an
arrow from $\D$ to $\D'$ if 
$\d(\D)=\d(\D')+1$, and $\D\lrar\D'$.

Let $\Q$ be a quiver. One can consider the {\sl quiver algebra}
$\F\Q$ (see e.g. [Vyb1]).

\subhead {1.2.9} \endsubhead
We are especially interested in the following quotients
of the quiver algebra $\F\Q(X,\d)$.

\definition{Definition A} 
The algebra $A(X,\d)$ is the quotient of 
$\F \op{Q}(X,\d)$ by the ``chain complex'' relations:

If $\D'$ is any simplex such that
$\d(\D')=k+1$ and $\D''$ is any simplex such that 
$\d(\D'')=k-1$, then
$$\sum\limits_{\matrix
\D : \d(\D)=k,\\
\D'\lrar\D\lrar\D''\endmatrix}
a(\D ,\D'') a(\D',\D )=0,$$
where $a(\D ,\D'')$, $a(\D',\D )$ are the generators of $A(X,\d)$.
\enddefinition

\definition{Definition B}
The algebra $B(X,\d)$ is the quotient of
$\F \op{Q}(X,\d)$ by the ``equivalence'' relations:

Suppose that $\D'$, $\D''$, $\D_1$ and $\D_2$ are simplices
of $X$  such that:
\roster
\item $\d(\D')=k+1$, $\d (\D'')=k-1$ and $\d(\D_1)=\d(\D_2)=k$, 
\item  $\D'\lrar\D_1\lrar\D''$ and $\D'\lrar\D_2\lrar\D''$.
\endroster
Then 
$$
b(\D_1, \D'') b(\D',\D_1)=
b(\D_2, \D'') b(\D',\D_2),
$$
where $b(\D_1, \D'')$, $b(\D',\D_1)$,
$b(\D_2, \D'')$, and $b(\D',\D_2)$
are the generators of $B(X,\d)$.
\enddefinition

\subhead {1.2.10} \endsubhead
We will present now an algebraic  interpretation of the category 
$\PP(X,\d)$ and the category $\RR(X,\d)$. 

Let us construct 
a functor $\Xi_A(X,\d):\PP(X,\d)\to\moda(X,\d)$
(resp. $\Xi_B(X,\d):\RR(X,\d)\to\modb(X,\d)$). Let 
$\S\in\PP(X,\d)$ (resp. $\RR(X,\d)$). Let $M=\Xi_A(X,\d)(\S)$
(resp. $\Xi_B(X,\d)(\S)$).
As a vector space $M=\oplus S(\D)$, where the sum is taken over all
simplices of $X$.  The action of the respective
algebras is defined as follows:
\roster
\item $e(\D)m=m$ for $m\in S(\D)$, and $e(\D)m=0$ for $m\in S(\D')$, 
$\D\neq\D'$,
\item $a(\D,\D')m=s(\D,\D')(m)$ for $m\in S(\D)$,
and $a(\D,\D')m=0$ for $m\in S(\D'')$, $\D\neq\D''$ 
(resp. $b(\D,\D')m=s(\D,\D')(m)$ for $m\in S(\D)$,
and $b(\D,\D')m=0$ for $m\in S(\D'')$, $\D\neq\D''$).
\endroster
 
\proclaim{Lemma} 
\roster
\item The functor
$$
\Xi_A(X,\d):\PP(X,\d)\to\moda(X,\d)
$$
is an isomorphism of categories.
\item The functor
$$
\Xi_B(X,\d):\RR(X,\d)\to\modb(X,\d)
$$
is an isomorphism of categories.
\endroster
\endproclaim

\demo{Proof} The proof is left to the reader.
\enddemo

\subhead {1.2.11} \endsubhead
If $\d=\text{bottom perversity}$, 
then the symbol $\d$ will just be omitted from the notation,
for example $\L(X,\d)$ 
will be denoted by $\L(X)$.

\subhead 1.3. Elements of Sheaf Theory\endsubhead
In this section we list some sheaf theoretic
definitions and constructions mostly lifted
from [BBD, Iver, KS, Pol], where we refer the
reader for many more details. 

\subhead 1.3.1 \endsubhead
Let $(X,\TT)$ be a topological space.
The set $\TT$ has a natural partial order:
for $U,V\in\TT$ we say that $U\leq V$
if $U\subeq V$. Let $\S$ be a presheaf on $\TT$,
and let $U\subeq V$. The restriction map
$S(V)\to S(U)$ is sometimes denoted by $r^V_U$,
and if $s\in S(V)$, then $r^V_U(s)$ is sometimes
denoted by $s|_U$.

We record the following standard definition here only to be able to
refer to axioms S1 and S2 later in the text.

\definition{Definition} A presheaf $\S$ on $\TT$ is called a 
{\sl sheaf} if it satisfies conditions S1 and S2 below.

S1. For any open set $U\subeq X$, any open covering $U=\cup_\a U_\a$,
any section $s\in S(U)$, $s|_{U_{\a}}=0$ for all $\a$ implies $s=0$.

S2. For any open set $U\subeq X$, any open covering $U=\cup_\a U_\a$,
any family $s_\a\in S(U_\a)$
satisfying $s_{\a}|_{U_\a\cap U_\b}=
s_{\b}|_{U_\a\cap U_\b}$ 
for all pairs $(\a,\b)$, there exists $s\in S(U)$
such that $s|_{U_{\a}}=s_{\a}$ for all $\a$.
\enddefinition 

The category of sheaves $\SH(\TT)$ on $(X,\TT)$ is a full 
subcategory in the category of presheaves $\Pre\TT$.

\subhead{1.3.2}\endsubhead 
Let $\S$ be a presheaf on $\TT$. 
Let  $S_x$ be the stalk of $\S$ 
at a point $x\in X$.
The image of $s\in S(U)$ in $S_x$ is denoted by $s_x$.
For $s\in S(U)$ its support is denoted by $\op{supp}(s)$.

\proclaim{1.3.3. Lemma [KS, 2.2.3]} Given a presheaf $\S$ on 
$(X,\TT)$,
there exists a sheaf $\S^+$ and a morphism $\theta:\S\to\S^+$.
Moreover, $(\S^+,\theta)$ is unique up to an isomorphism, and
for any $x\in X$, $\theta_x: S_x\to S_x^+$ is an isomorphism.
\endproclaim

Sometimes $\S^+$ is called the sheafification of $\S$.
The sheafification functor is denoted by $+:\Pre\TT\to\SH(\TT)$.

\subhead{1.3.4}\endsubhead 
Let $A$ be a vector space.
The constant sheaf assosiated to $A$ is denoted 
by $\A$. 

\subhead 1.3.5 [KS, 2.3.1, 2.5], [Iver, II.6.6]\endsubhead 
Let $X$ and $Y$ be two topological spaces, $f:Y\to X$ a continuous map.
Let $\T$ be a sheaf on $Y$. The 
direct image of $\T$ by $f$ (resp. direct image with proper supports)
is denoted by $f_*\T$ (resp. $f_!\T$).
Let $\S$ be a sheaf on $X$. The inverse image of $\S$ by $f$
(resp. inverse image with proper supports) is
denoted by $f^*\S$ (resp. $f^!\S$).
Let $f:Y\hook X$ be an inclusion of a subset and let $\S$ be 
a sheaf on $X$. Then $f^*\S$ is also denoted by $\S|_Y$.

\definition{1.3.6. Definition}
Let $X$ be a finite connected simplicial complex
equipped with the usual topology (cf. 1.1),
and let
$\D\subeq X$ be a simplex in $X$.
A sheaf $\S$ of $\F$-vector spaces on $X$ is called 
{\sl constructible with respect to
the triangulation}, if $\S|_{\D}$
is the constant sheaf associated to a finite dimensional 
$\F$-vector space for all $\D$.
\enddefinition

The category of constructible sheaves (denoted by $\SH_c(X)$) is a 
full subcategory of $\SH(\TT)$, where $\TT$ is the usual topology on 
$X$.

\subhead{1.3.7}\endsubhead
Let $\Cal A$ be an abelian category. The bounded homotopic 
(resp. derived) category of $\Cal A$ will be denoted by
$\CC(\Cal A)$ (resp. $\DD(\Cal A)$). 

\definition{Definition} Let $X$ be
a finite connected simplicial complex equipped with the
usual topology $\TT$. Let $(\S^i,d^i)\in\DD(\SH(\TT))$.
Then $(\S^i,d^i)$ is called {\sl cohomologically constructible}
with respect to the triangulation
if the sheaves 
$\bold H^i(\S^\bu)=\bold{Ker}d^i/\bold{Im}d^{i-1}$
are constructible with respect to the trangulation.
The category of cohomologically constructible complexes of 
sheaves (denoted by $\DD_c(X)$) is a 
full subcategory of $\DD(\SH(\TT))$.
\enddefinition

The following theorem is due to M. Kashiwara [Kash] and A. Shepard
[Shep]. We refer the reader to [KS, 8.1.11].

\proclaim{Theorem}
The following categories
$$
\DD(\SH_c(X))\simeq\DD_c(X)
$$
are equivalent.
\endproclaim 

\subhead{1.3.8 [BBD, 1.3.1]}\endsubhead 
Let $\Cal D$
be a triangulated category, and let 
$(\Cal D^{\leqslant 0},\Cal D^{\geqslant 0})$
be a $t$-structure on $\Cal D$ [BBD, 1.3.1].
The {\sl core} of the $t$-structure 
$(\Cal D^{\leqslant 0},\Cal D^{\geqslant 0})$
is the full subcategory 
$\Cal A:=\Cal D^{\leqslant 0}\cap\Cal D^{\geqslant 0}$.
It is known [BBD, 1.3.6] that the core of a $t$-structure
is an abelian category stable under extensions.

\subhead{1.3.9}\endsubhead
This subsection is mostly lifted from [BBD].
Let $X$ be a topological space, and $X=\sqcup_{S\in \Cal S} S$ 
be a finite decomposition (stratification)
of $X$ into the disjoint union of locally closed subspaces (strata).
Let $p$ be a BBDG perversity. Let $\DD(X)$ be the derived category
of (cohomologically constructible with respect to some triangulation
refining the stratification) sheaves of $\F$-vector spaces on $X$.

\definition{Definition [BBD, 2.1.2]}
The subcategory ${}^{p}\Cal D^{\leqslant 0}(X)$
(resp. ${}^p\Cal D^{\geqslant 0}(X)$) of $\DD(X)$
is the subcategory formed by the complexes $\S^\bu$
such that for every inclusion of a stratum
$i_S:S\hook X$, we have 
$\bold {H}^n(i^*_S \S^\bu)=0$ for $n>p(\dim S)$
(resp. $\bold {H}^n(i^!_S \S^\bu)=0$ for $n<p(\dim S)$).
\enddefinition

\proclaim{Proposition [BBD, 2.1.4]}
$({}^p\Cal D^{\leqslant 0}(X),{}^p\Cal D^{\geqslant 0}(X))$
is a $t$-structure on $\DD(X)$.
\endproclaim

The core of this $t$-structure, denoted by
$\Cal M_{\Cal S}(X,p):=
{}^p\Cal D^{\leqslant 0}(X)\cap{}^p\Cal D^{\geqslant 0}(X)$
is called the category of 
{\sl cohomologically constructible with respect to the 
stratification $p$-perverse sheaves} on $X$.

In particular, if the stratification $\Cal S$ is a triangulation 
$\TT$ of $X$,
we sometimes denote the category $\Cal M_{\TT}(X,p)$ of 
the $p$-perverse sheaves  cohomologically constructible with 
respect to the triangulation simply by $\M(X,p)$. 

Let $\TT$ be a triangulation refining a stratification $\Cal S$.
Then there is an exact embedding functor 
$\refn:\Cal M_{\Cal S}(X,p)\hookrightarrow\Cal M_{\TT}(X,p)$
(see [BBD, 2.1.14, 2.1.15]), i.e.
$\Cal M_{\Cal S}(X,p)$ may be considered as a subcategory of 
$\Cal M_{\TT}(X,p)$.    

Let $\d$ be the cellular perversity corresponding to $p$ (see 1.1).
Below we will often denote the category $\Cal M_{\SS}(X,p)$ by 
$\Cal M_{\SS}(X,\d)$.

\subhead{1.3.10}\endsubhead
It is known that the simple objects of $\M(X,\d)$ are in
one-to-one correspondence with simplices $\D$ in $X$.
In this subsection we will give an explicit description 
of the simple objects $\S_\D$ in 
$\M(X,\d)=\M(X,p)$, following [Pol].

\proclaim{Lemma [Pol, 1.1]} For a simplex $\D$ in $X$ we have:
$$
\S_\D=\cases
i_{\D*}\FF_\D[-p(\dim\D)], & \dim\D \text{ is of type } * \\
i_{\D!}\bold{or}_\D[-p(\dim\D)], & \dim\D \text{ is of type } !\ , \\
\endcases
$$
where $\FF_\D$ is the constant sheaf on $\D$ associated to
$\F$, and $\bold{or}_\D$ is the orientation sheaf on $\D$
(see [KS, 3.3]).
\endproclaim

\head{Chapter 2.  
Constructible sheaves and perverse topology}
\endhead

\subhead 2.1. Sheaves and a base for topology \endsubhead
In this section we discuss the relationship between 
sheaves on a topological space $(X,\TT)$ and presheaves on a 
base for topology. Many technical details are lifted from [KS]. 

\subhead {2.1.1. [Jos]}\endsubhead
Let $(X,\TT)$ be a topological space. A subfamily $\CB$ of $\TT$
is said to be a base for $\TT$ if every member of $\TT$ can be
expressed as the union of some members of $\CB$.

\proclaim{Lemma} Let $X$ be a set and $\CB$ a family of its subsets
covering $X$. Then the following statements are equivalent:
\roster
\item there exists a topology on $X$ with $\CB$ as a base,
\item for any $B_1,B_2\in\CB$, $B_1\cap B_2$ can be expressed as the
union of some members of $\CB$,
\item for any $B_1,B_2\in\CB$ and $x\in B_1\cap B_2$, there exists
$B_3\in\CB$ such that $x\in B_3$ and $B_3\sub B_1\cap B_2$.
\endroster
\endproclaim

\subhead {2.1.2}\endsubhead
Let $(X, \TT)$ be a topological space and let $\CB$ be a base for
$\TT$. The set $\CB$ has a natural partial order: for $U,V\in\CB$ 
we say that $U\leq V$ if $U\subeq V$. 
Let $W\in\TT$ be an open set.
Define 
$$
\CB|_W=\{ U\sub W : U\in\CB \}.
$$
It is sometimes convenient to assign some indices to the 
elements of 
$\CB|_W$, 
i.e. we regard $\CB|_W=\{U_\l\}$ as a family of subsets 
parametrized by $\l$.

\subhead {2.1.3}\endsubhead
In this subsection we construct a functor{\ }  
$\wti{\ } : \Pre\CB\to\Pre\TT$.

On objects: let $\S\in\Pre\CB$ be a presheaf on $\CB$.
Let $W\in\TT$ be an open set and let $\CB|_W=\{U_\l\}$. 
We define
$$
\ti S(W)=\vplim_{\l} S(U_\l)=
\{s\in\prod_\l S(U_\l) : s(\l)|_{U_\m}=s(\m) \text{ for }
U_\m\subeq U_\l \}.
$$
Sometimes we write $s_\l$ instead of  $s(\l)$. 

We have a natural restriction map $\ti S(W)\to S(U_\l)$:
for $s\in\ti S(W)$, 
$r^W_{U_\l}(s)=s(\l)$.

Let $W_1,W_2\in\TT$ be such open sets that $W_2\subeq W_1$,
and let $\CB|_{W_1}=\{U_\l\}$ and
$\CB|_{W_2}=\{U_\m\}$.  For $s\in\ti S(W_1)$ we define:
$$
r^{W_1}_{W_2}(s)=\prod_\m r^{W_1}_{U_\m}(s)=\prod_\m s(\m).
$$
We have $s(\m)|_{U_\n}=s(\n)$ for $U_\n\subeq U_\m$ since
$s\in\vplim_{\l} S(U_\l)$. Therefore 
$r^{W_1}_{W_2}(s)\in\vplim_{\m} S(U_\m)=\ti S(W_2)$.
Clearly $r^{W}_{W}=\Id$.

Let $W_1, W_2, W_3\in\TT$ be such open sets that 
$W_3\subeq W_2\subeq W_1$. Let $\CB|_{W_1}=\{U_\l\}$,
$\CB|_{W_2}=\{U_\m\}$, and $\CB|_{W_3}=\{U_\n\}$.
Let $s\in\ti S(W_1)=\vplim_{\l} S(U_\l)$. Then
$r^{W_1}_{W_3}(s)=\prod_\n s(\n)$ and
$r^{W_2}_{W_3}\circ r^{W_1}_{W_2}(s)=
r^{W_2}_{W_3}(\prod_\m s(\m))=\prod_\n s(\n)$. Thus
$$
r^{W_1}_{W_3}=r^{W_2}_{W_3}\circ r^{W_1}_{W_2}.
$$
We leave it to the reader to construct the functor on morphisms.

\subhead {2.1.4} \endsubhead
Let $|_{\CB} :\Pre\TT\to\Pre\CB$ be the obvious ``restriction'' 
functor.

\proclaim{Lemma} $|_{\CB}\circ\wti{\ }\simeq\Id$.
\endproclaim

The proof is left to the reader.

\subhead {2.1.5} \endsubhead
Let{\ } $+ : \Pre\TT\to\SH(\TT)$ be the
sheafification functor and $\imath:\SH(\TT)\to\Pre\TT$ be the
inclusion functor. We define $\php: \Pre\CB\to\SH(\TT)$, 
$\php=+\circ\wti{\ }$, and $\phm:\SH(\TT)\to\Pre\CB$,
$\phm=|_{\CB}\circ\imath$.

\proclaim{2.1.6. Lemma} $\php\circ\phm\simeq\Id$.
\endproclaim

The proof is left to the reader.

\subhead 2.2. Proof of the Theorem A\endsubhead
In this section we give the sequence of elementary lemmas
leading to the proof of Theorem A. This theorem
is well known [Kash, Shep, KS]. The proofs of the lemmas are left 
to the reader.

\subhead {2.2.1}\endsubhead
Let $X$ be a finite connected simplicial complex.
Let us consider the partially ordered set $\L(X,\d)$ 
introduced in 1.2, when $\d=\text{bottom perversity}$.
In this case we will denote it by $\L(X)$. Clearly $\D\geq\D'$
if and only if $\D\subeq\D'$. For every simplex $\D$ we denote
its open star by $\D^*$. 

\subhead {2.2.2}\endsubhead
This subsection is mostly lifted from [KS, 8.1.4]. 
Let $K$ be a finite simplicial set with a set of vertices $\ver$ 
(cf. 1.1). Let $\D$ be a simplex of $|K|$. Let $x\in\D$ be a point.
For $0<\e\leq 1$, set $I_\e=\{\a\in\R : \e\leq \a\leq 1\}$
and define the map $\pi_\e: I_\e\times\D^*\to\D^*$ by:
$$
\pi_\e(\a,y)(v)=\a y(v)+(1-\a)x(v)
\qquad
\text{ for } v\in\ver.
$$
The map $\pi_\e$ is continuous and surjective, $\pi_\e(1,\cdot)$
is the identity $\{1\}\times\D^*\simeq\D^*$, and $\pi_\e(\e,\cdot)$
is a homeomorphism $\{\e\}\times\D^*\overset{\sim}\to\to
\pi_\e(\{\e\}\times\D^*)$. We denote 
$\D^*(x,\e):=\pi_\e(\{\e\}\times\D^*)$. In other words:
$$ 
\D^*(x,\e)=\{\e y+(1-\e)x : y\in\D^*\}
$$

Clearly, $\D^*(x,\e)\sub\D^*$ for $0<\e<1$. A set $\D^*(x,\e)$
will be called the {\sl $\e$-neighborhood} of a point $x\in\D$.
It is clear that the $\e$-neighborhoods form a base for topology of 
$|K|$.

Let $X$ be a finite connected simplicial complex. For a point 
$x\in\D$, (and $0<\e\leq 1$) its $\e$-neighborhood $\D^*(x,\e)$ 
is the homeomorphic
image of the corresponding $\e$-neighborhood in $|K|$.
Certainly the $\e$-neighborhoods form a base $\CB$ for topology 
$\TT$ of $X$.
 
\subhead {2.2.3}\endsubhead
Let $\T$ be a presheaf on $\L(X)$. Let $\D^*(x,\e)$ be some 
$\e$-neighborhood of a point $x\in\D$. Notice that if $\D'$ 
is another simplex, then $\D'\cap\D^*(x,\e)\neq\ves$ if and only
if $\D'\leq\D$.

\definition{Definition} (a) A function 
$\vp:\D^*(x,\e)\to\oplus_\D T(\D)$ is called 
{\sl $\T$-constructible} if:
\roster
\item $\vp|_{\D'\cap\D^*(x,\e)}\equiv\vp(\D')\in T(\D')$,
is a constant function for $\D'\leq\D$,
\item $\vp(\D'')=t(\D',\D'')\vp(\D')$  for  
$\D''\le\D'\leq\D$.
\endroster

(b) Let $U\subeq X$ be an open set. A function 
$\vp: U\to\oplus_\D T(\D)$ is called {\sl locally 
$\T$-constructible} if every point $x\in U$ has an
$\e$-neighborhood such that $\vp$ becomes $\T$-constructible
when restricted to this neighborhood. Notice that 
a locally $\T$-constructible function becomes constructible
upon restriction to any $\e$-neighborhood of any point.
\enddefinition

\proclaim{Lemma} Let $x\in\D\sub X$ be a point and let 
$0<\e\leq 1$. Then
$$
\{\text{the space of }\T
\text{-constructible functions on }\D^*(x,\e)\}=T(\D).
$$
Moreover, if $\D_1^*(x,\e_1)\supeq\D_2^*(y,\e_2)$ then    
the map from the 
space of $\T$-constructible functions on $\D_1^*(x,\e_1)$
to the 
space of $\T$-constructible functions on $\D_2^*(y,\e_2)$
given by restriction is precisely $t(\D_1,\D_2)$.
\endproclaim

\subhead {2.2.4}\endsubhead  
Let $\T$ be a presheaf on $\L(X)$. 
We will construct a presheaf (denoted by the same symbol $\T$)
on $\CB$ assigning to each $\D^*(x,\e)$ the space of 
$\T$-constructible functions on it, with restriction maps 
given by restricting functions. Such a presheaf will be called 
constructible. (In other words, $\T\in\Pre\CB$ is constructible
if $\Id=t(\D^*,\D^*(x,\e)):T(\D^*)\overset\sim\to\to T(\D^*(x,\e))$ 
for all appropriate $\e$ and $x\in\D$.)
The full subcategory of $\Pre\CB$ consisting
of constructible presheaves will be denoted $\Pre_c\CB$.
Clearly, $\Pre\L(X)=\Pre_c\CB$. This explains the same notation for
the objects of the two categories. Sometimes we will use one
instead of the other.

\proclaim{2.2.5. Lemma} Let $\T$ be a presheaf on $\L(X)$.
Let $\ti\T$ be a presheaf on $\TT$ corresponding to 
$\T\in\Pre_c\CB$
(see 2.1.3). Let $U\in\TT$ be an open set. Then    
$$
\ti T(U)=\{\text{the space of locally }\T
\text{-constructible functions on } U\}.
$$
Moreover, the restriction map $\ti T(U_1)\to\ti T(U_2)$,
for $U_2\subeq U_1$ applied to a function
$\vp\in T(U_1)$ is the restriction of $\vp$
to $U_2$.
\endproclaim

\proclaim{2.2.6. Lemma}(a) $\ti\T$ is a sheaf.

(b) $\ti\T$ is constructible with respect to the triangulation.
\endproclaim

\subhead 2.2.7 \endsubhead
In this subsection we will slightly generalize the argument
of 2.2.6(b).  
Let $Y\subeq X$ be a closed 
union of simplices. Let $\T$ be a presheaf on $\L(X)$
and let $\ti\T$ be the corresponding constructible sheaf.
We can restrict the presheaf $\T$ to a presheaf $\T|_Y$ 
on the partially ordered set $\L(Y)$ in the obvious way.

\proclaim{Lemma} $\ti\T|_Y=\wti{\T|_Y}$.
\endproclaim

Let $U\subeq X$ be an open union of simplices. It is clear that
$\ti\T|_U=\wti{\T|_U}$. (Here $\wti{\T|_U}$ is the obvious 
subsheaf of $\wti{\T|_{\ov U}}$.) 
Then the statement of the Lemma holds for any 
locally closed $Y\subeq X$.

\subhead {2.2.8} \endsubhead
Lemma 2.2.6 implies that the functor{\ } $\wti{\ }$ (see 2.1)
restricts to a functor $\psp: \Pre\L(X)=\Pre_c\CB\to\SH_c(X)$.
Notice that if $\psp$ is considered as a functor from 
$\Pre_c\CB$, then it coincides with the restriction of the
functor $\php$ defined in 2.1.
The functor $\psm:\SH_c(X)\to\Pre\L(X)$ is constructed as follows.
Let $\S\in\SH_c(X)$ and let $\T=\psm(\S)$. Then
\roster
\item $T(\D)=S(\D^*)$,
\item $t(\D_1,\D_2)=r^{\D_1^*}_{\D_2^*}$ (restriction maps of $\S$).
\endroster
It is clear from our construction that $\psm\circ\psp\simeq\Id$.
(Alternatively, we can use Lemma 2.1.4.)

\proclaim{2.2.9. Lemma} $\psp\circ\psm\simeq\Id$.
\endproclaim

\subhead {2.2.10} \endsubhead
Summarizing, we have constructed the functor $\psp:\Pre\L(X)\to\SH_c(X)$
which to a presheaf $\T$ assignes the sheaf of locally
$\T$-constructible functions, and the functor
$\psm:\SH_c(X)\to\Pre\L(X)$ which to a constructible sheaf $\S$ 
assignes the presheaf of sections on open stars of simplices.
We have proved that these functors are 
quasi-inverse to each other. In 1.2 we constructed 
the equivalence functors:
$$
\Pre\L(X)@>{F}>>\RR(X)@>{\Xi_B(X)}>>\modb(X).
$$ 
Let $\Th(X):\SH_c(X)\to\modb(X)$
be the equivalence functor defined as 
$\Th(X)=\Xi_B(X)\circ F\circ\psm$. 
In other words, we have

\proclaim{Theorem A} The functor
$$
\Th(X):\SH_c(X)@>{\sim}>>\modb(X)   
$$
is an equivalence of categories.
\endproclaim

\subhead 2.3. Perverse topology
\endsubhead
In this section we introduce perverse topology and
prove Theorem B1. 
For notational conventions on simplices, simplicial sets,
and simplicial complexes see 1.1.

\subhead 2.3.1 \endsubhead
Let $K$ be a finite simplicial set. Let $\D\in K$ be
a simplex and $|\D|$ be the corresponding geometric simplex. 
Let $|\D|^*$ be the open star of $|\D|$ (i.e. the union of all simplices
having $|\D|$ as a face). 
If $\D_1$ and $\D_2$ are two simplices, then
$$
|\D_1|^*\cap|\D_2|^*=
\cases
|\D_1\cup\D_2|^*, &
\text{ if } \D_1\cup\D_2\in K \\
\emptyset, &\text{ otherwise. } \\
\endcases
$$

\subhead 2.3.2\endsubhead 
Let $\D\in K$ be a simplex and let $\ov{|\D|}$
be the closure of the corresponding geometric simplex.
Notice that
$$
\ov{|\D|}=\bigsqcup_{\D'\subeq\D}|\D'|.
$$
If $\D_1$ and $\D_2$
are two simplices, then
$\ov{|\D_1|}\cap\ov{|\D_2|}=\ov{|\D_1\cap\D_2|}$.
(We assume that $\ov{|\emptyset|}=\emptyset$.)

\subhead 2.3.3 \endsubhead
Until the end of this section 
$X$ is assumed to be a finite connected simplicial complex.
Let us fix a perversity $\d$.
Let $\L(X,\d)$ be the partially ordered set introduced in 1.2.2.
For a simplex $\D\sub X$ we define
the {\sl perverse} star $U(\D,\d)$ of $\D$ to be the 
following union of simplices:
$$
U(\D,\d)=\bigsqcup_{\D'\leq\D}\D'\ .
$$
If $\d=\text{bottom perversity}$, then $U(\D,\d)=\D^*$,
and if $\d=\text{top perversity}$, then $U(\D,\d)=\ov\D$.
Let $\CB(X,\d)$ be the family of sets $U(\D,\d)$
parametrized by all simplices of $X$.

\proclaim{Lemma} There exists a topology on $X$ with 
$\CB(X,\d)$ as a base.
\endproclaim
\demo{Proof} First of all the family of subsets $\CB(X,\d)$
covers $X$. Indeed, if $x\in X$ is a point, then $x\in\D$
for some $\D$, and  $x\in U(\D,\d)$.

Let $\D$ and $\D'$ be two simplices. Then 
$$
U(\D,\d)\cap U(\D',\d)=\bigcup_{\D''\leq\D,~~\D''\leq\D'}
U(\D'',\d).
$$
Indeed, 
if $\D''\subeq U(\D,\d)\cap U(\D',\d)$ 
(i.e. $\D''\leq\D$ and $\D''\leq\D'$),
then 
$\D''\subeq U(\D'',\d)$. 
Vice versa
if $\D''\leq\D$ and $\D''\leq\D'$,
then $U(\D'',\d)\subeq U(\D,\d)\cap U(\D',\d)$ since
if $\ti\D\subeq U(\D'',\d)$, then $\ti\D\leq\D''\leq\D$,
and $\ti\D\leq\D''\leq\D'$.

Thus, the claim follows from Lemma 2.1.1.
\enddemo

The topology generated by $\CB(X,\d)$ will be denoted by 
$\TT(X,\d)$. An element $U\in\TT(X,\d)$ will be called a 
{\sl perverse set}.

\subhead 2.3.4 \endsubhead
The partially ordered set $\L(X,\d)$ can be interpreted 
now as a set of perverse stars $U(\D,\d)$ with 
the partial order given by inclusion. (Indeed, 
$U(\D',\d)\subeq U(\D,\d)$ if and only if $\D'\leq\D$.)
In other words, $\L(X,\d)=\CB(X,\d)$ as partially ordered sets
and $\Pre\L(X,\d)=\Pre\CB(X,\d)$. Let $\S\in\Pre\CB(X,\d)$
and let $\ti\S$ be a presheaf on $\TT(X,\d)$ associated to $\S$ 
(see 2.1).

\proclaim{Lemma} $\ti\S$ is a sheaf.
\endproclaim

\demo{Proof} Let $U$ be a perverse set, and  let 
$U=\cup_{\a=1}^l U_\a$ be a covering of $U$ by perverse sets. 
We claim that each $U(\D,\d)\subeq U$ is contained in some $U_\a$.
Indeed, if $U(\D,\d)\subeq\cup_\a U_\a=\cup_\l U(\D_\l,\d)$, then
$\D\subeq\cup_\l U(\D_\l,\d)$. Then there exists such $\l_0$ that
$\D\subeq U(\D_{\l_0},\d)$. Then 
$U(\D,\d)\subeq U(\D_{\l_0},\d)\subeq U_\a$
for some $\a$.

We have to check axioms 1.3.1.S1 and 1.3.1.S2.

1.3.1.S1. Let $s\in\ti S(U)$. This means that 
$s=\prod_{\D\subeq U} s(\D)$ with $s(\D)\in S(U(\D,\d))$
and $s(\D)|_{U(\D',\d)}=s(\D')$
for $U(\D',\d)\subeq U(\D,\d)$. Let $s|_{U_\a}=0$
for all $\a$. This implies that $s(\D)=0$ whenever
$U(\D,\d)\subeq U_\a$ for some $\a$, but since every perverse star
is a subset of some $U_\a$, $s(\D)=0$ for all 
$\D\subeq U$ and $s=0$.

1.3.1.S2. Let $s_\a\in\ti S(U_\a)$, and let
$s_\a|_{U_\a\cap U_\b}=s_\b|_{U_\a\cap U_\b}$ for all $\a,\b$.
Each $s_\a=\prod_{\D\subeq U_\a} s_\a(\D)$ with 
$s_\a(\D)\in S(U(\D,\d))$.
Let $U(\D,\d)\subeq U$. Then there exists such $\a$ that
$U(\D,\d)\subeq U_\a$. We define 
$s(\D)=s_\a|_{U(\D,\d)}=s_\a(\D)$.
This definition does not depend on $\a$ due to compatibility
on intersections. We set $s=\prod_{\D\subeq U} s(\D)$. 
Clearly, $s|_{U_\a}=s_\a$.
\enddemo

\subhead 2.3.5 \endsubhead
Recall that in 2.1 we constructed functors{\ } 
$\wti{\ }: \Pre\CB(X,\d)\to \Pre\TT(X,\d)$,
$\php:\Pre\CB(X,\d)\to\SH(\TT(X,\d))$, and
$\phm:\SH(\TT(X,\d))\to\Pre\CB(X,\d)$. 
Lemma 2.3.4 implies that in the case of perverse topology
the functors{\ } $\wti{\ }=\php$ are isomorphic.

\proclaim{Lemma} The functors $\php$ and $\phm$ 
are quasi-inverse to each other. 
\endproclaim

\demo{Proof} Follows from Lemmas 2.1.4, 2.1.6.
\enddemo

\subhead 2.3.6 \endsubhead
Recall that in 1.2 we constructed the equivalence functors:
$$
\Pre\L(X,\d)@>{F}>>\RR(X,\d)@>{\Xi_B(X,\d)}>>\modb(X,\d).
$$ 
Let $\Omega(X,\d):\SH(\TT(X,\d))\to\modb(X,\d)$
be the equivalence functor defined as follows:
$\Omega(X,\d)=\Xi_B(X,\d)\circ F\circ\phm$. 
Summarizing the discussion we formulate: 

\proclaim{Theorem B1}The functor
$$
\Omega(X,\d):\SH(\TT(X,\d))@>{\sim}>>\modb(X,\d)   
$$
is an equivalence of categories.  
\endproclaim

\subhead 2.4. Borel-Moore-Verdier duality for simplicial sheaves and 
cosheaves
\endsubhead
In this section we list some basic facts about simplicial sheaves 
and cosheaves. The material is essentially borrowed
from [Shep] and [GMMV], where we refer the reader for 
many more details.  

In this section $X$ is assumed be a finite connected simplicial 
complex.
We denote by $V^*=\Hom_\F(V,\F)$, where $V$ is a vector space.
If $f$ is a linear map, then $f^*$ denotes the adjoint map.

\subhead{2.4.1}\endsubhead
Let $\S_1$ and $\S_2$ be two simplicial sheaves. Then
the sheaf $\T=\HOM(\S_1,\S_2)$ is defined as follows:
$$
T(\D)=\Hom(\S_1|_{\D^*},\S_2|_{\D^*}).
$$
If $\D_1\subeq\D_2$, then the map $t(\D_1,\D_2):T(\D_1)\to T(\D_2)$
is given by restriction of a map 
$\S_1|_{\D_1^*}\to\S_2|_{\D_1^*}$ to a map
$\S_1|_{\D_2^*}\to\S_2|_{\D_2^*}$ (note that $\D_2^*\subeq\D_1^*$).

\definition{2.4.2. Definition} Let $\D\sub X$ be a simplex
and let $V$ be a vector space. We define the simplicial sheaf
$[\D]^V$ as follows:
$$
[\D]^V(\D')=\cases
V, & \D'\subeq\D \\
0, & \D'\not\subeq\D, \\
\endcases
$$
where the restriction maps between copies of $V$ are all
identity maps. We will also write $[\D]=[\D]^\F$.
It is easy to see that for any $\D$ and $V$,
 $[\D]^V$ is an injective simplicial sheaf. 
\enddefinition

\proclaim{Lemma} Let $\S$ be a simplicial sheaf, and let
$[\D]^V$ be an injective sheaf. Then
\roster
\item $\Hom(\S,[\D]^V)=\Hom_\F (S(\D),V)$.
\item $\HOM(\S,[\D]^V)=[\D]^{\Hom_\F(S(\D),V)}$.
\item In particular, if $\D''\subeq\D'$, then
$\Hom([\D']^V,[\D'']^W)=\Hom_\F(V,W)$.
\endroster
\endproclaim

\demo{Proof} The proof is left to the reader.
\enddemo

\subhead{2.4.3} \endsubhead
Let $\S$ be a simplicial sheaf. There
is a canonical map $\S\to [\D]^{S(\D)}$ induced by
$\Id:S(\D)\to S(\D)$. Taking the sum over all simplices we obtain 
a canonical map $\S\to\oplus_{\D}[\D]^{S(\D)}=\I^0$ which
is clearly injective. Moreover,

\proclaim{Lemma} For any simplicial sheaf $\S$ there exists
the canonical injective resolution, $\S@>{q.i.}>>\I^\bullet$.
\endproclaim

\demo{Proof} See [Shep, 1.3]. \enddemo

\remark{Remark} Let $\FF$ be a constant sheaf on $X$.
The canonical injective resolution of $\FF$ is as follows:
$$
\FF@>>>\oplus_{\D_0}[\D_0]
@>>>\oplus_{\D_1\sub\D_0}[\D_1]
@>>>\oplus_{\D_2\sub\D_1\sub\D_0}[\D_2]
@>>>\dots
$$
Taking the global sections (see [Shep, 1.4]) we obtain the complex
of vector spaces computing the simplicial cohomology
(with coefficients in $\F$) of the first barycentric 
subdivision $\hx$.
\endremark

\subhead{2.4.4}\endsubhead
Let $\CC(\Pre\L(X))$ 
(resp. $\DD(\Pre\L(X))$) be 
the homotopic (resp. derived) category 
of the category of simplicial sheaves.

\definition{Definition} The dualizing complex
$\om_X=\om_X^\bullet\in \DD(\Pre\L(X))$ has in degree $-i$
the sheaf: 
$$
\om_X^{-i}=\bigoplus_{\dim\D=i}[\D].
$$ 
The boundary map $\om_X^{-i}\to\om_X^{-i+1}$ is the
zero map between components $[\D]$ and $[\D']$
if $\D'\not\sub\D$, and is induced by multiplication by 
$[\D:\D']$ if $\D'$ is a $\op{codim}\ 1$ face of $\D$.
\enddefinition

\subhead{2.4.5}\endsubhead
In 1.2 we have seen that $\Pre\L(X)=\RR(X)$. Thus,
$\CC(\Pre\L(X))=\CC(\RR(X))$ and $\DD(\Pre\L(X))=\DD(\RR(X))$. 
We transport the injective objects $[\D]^V$ and the
dualizing complex $\om_X$ to the categories $\RR(X)$
and $\DD(\RR(X))$ respectively, preserving the
notation.

Let $(\S^i,d^i_S)\in\CC(\RR(X))$, 
and let us consider a bicomplex:
$$
D\S^{i,j}=\HOM(\S^{i},\om^{j}_X)=
\bigoplus_{\dim\D=-j}[\D]^{S^{i}(\D)^*}.
$$ 
The differential 
$$
\matrix\format\l \\
D\S^{i,j}@>{d_{I}^{i,j}}>> D\S^{i,j+1}, \\
\bigoplus_{\dim\D=-j}[\D]^{S^{i}(\D)^*}
@>{d_{I}^{i,j}}>>
\bigoplus_{\dim\D=-j-1}[\D]^{S^{i}(\D)^*} \\
\endmatrix
$$
is the zero map between components $[\D]^{S^{i}(\D)^*}$
and $[\D']^{S^{i}(\D')^*}$ if $\D'\not\sub\D$
and is induced by 
$[\D:\D']s^{i}(\D',\D)^*:S^{i}(\D)^*\to
S^{i}(\D')^*$, if $\D'$ is a $\codim\ 1$ face of
$\D$.

The differential 
$$
\matrix\format\l \\
D\S^{i,j}@>{d_{II}^{i,j}}>> D\S^{i-1,j}, \\
\bigoplus_{\dim\D=-j}[\D]^{S^{i}(\D)^*}
@>{d_{II}^{i,j}}>>
\bigoplus_{\dim\D=-j-1}[\D]^{S^{i-1}(\D)^*} \\
\endmatrix
$$
is the sum of the morphisms
$[\D]^{S^{i}(\D)^*}\to[\D]^{S^{i-1}(\D)^*}$ 
induced by the linear maps
$(-1)^{j-i+1}d_S^{i-1}(\D)^*:{S^{i}(\D)^*}\to
{S^{i-1}(\D)^*}$. 

Let us consider the complex
$$
D\S^p=\bigoplus_{p=j-i}D\S^{i,j}
$$ 
with the differential $d=d_{I}+d_{II}$.
The complex $D\S^\bu$ is called the Verdier dual of $\S^\bu$.
Notice that:
$$
D\S^\bu=\HOM^\bu(\S^\bu,\om^\bu_X),
$$
where we consider 
$\om^\bu_X$ as an object in $\CC(\RR(X))$.
Since $\om^\bu_X$ is a complex of injective objects,
the functor $D$ induces the functor 
$\DD(D):\DD(\RR(X))\to\DD(\RR(X))$.

\subhead {2.4.6}\endsubhead
Recall that $\RR(X)=\RR(X,\text{bottom perversity})$,
and let us denote $\RR'(X)=\RR(X,\text{top perversity})$.

In this subsection we construct the duality functor 
$*:\CC(\RR(X))\to\CC(\RR'(X))$
as follows. Let $(\S^\bu,d^\bu_S)\in\CC(\RR(X))$. We define:
\roster
\item $*S^i (\D)=S^{-i}(\D)^*$, 
\item if $\D'$ is a $\codim\ 1$ face of $\D$, then
$*s^{i}(\D,\D')=
s^{-i}(\D',\D)^*:S^{-i}(\D)^*\to
S^{-i}(\D')^*$,
\item $*d^i(\D)=(-1)^{i+1}
d_S^{-i-1}(\D)^*:S^{-i}(\D)^*\to S^{-i-1}(\D)^*$. 
\endroster
It is clear that 
the contravariant functor $*$ transforms quasi-isomorphisms to
quasi-isomorphisms,  and hence induces the functor
$\DD(*):\DD(\RR(X))\to\DD(\RR'(X))$. It is easy to see
that $\DD(*)$ is a contravariant equivalence of categories,
$\DD(\RR(X))^{\text{opp}}\simeq\DD(\RR'(X))$.

\subhead {2.4.7}\endsubhead
Recall that $\PP(X)=\PP(X,\text{bottom perversity})$,
and let us denote  $\PP'(X)=\PP(X,\text{top perversity})$.

\proclaim{Lemma} The following categories
$$
\split 
\RR(X)=\PP(X), \\
\RR'(X)=\PP'(X)
\endsplit
$$
are isomorphic.
\endproclaim

The proof is left to the reader.

The derived functor 
$\DD(\mp):\DD(\RR'(X))\to\DD(\PP'(X))$
is an equivalence of categories.

\head
{Chapter 3.
Sheaves constant along perverse simplices}
\endhead

\subhead 3.1. $\SH_c(X,\d)\simeq\SS(X,\d)$ \endsubhead
Let $X$ be a finite connected simplicial complex
and let $\hx$ be its first barycentric subdivision.
We fix a perversity $\d$.
For definitions and notation on perverse simplices
see 1.1.

\definition{3.1.1. Definition} A sheaf $\S$ on $X$ is called
constant along perverse simplices if for any perverse
simplex $Y$ the restriction $\S|_{Y}$ is a constant
sheaf associated to a finite dimensional $\F$-vector space.
\enddefinition
The category of sheaves constant along $(-\d)$-perverse simplices
is denoted by $\SH_c(X,\d)$. It is a full subcategory of 
the category $\SH_c(\hx)$

\definition{3.1.2. Definition} 
The category $\SS (X,\d)$
is a full subcategory of the abelian category $\Pre\L(\hx)$. 
An object $\T$ of $\Pre\L(\hx)$
belongs to $\SS(X,\d)$ 
if for any $\HD,\HD'\subeq\PD$ we have:
\roster
\item $T(\HD)=T(\HD')$,
\item if $\HD\subeq\HD'$, then
$t(\HD,\HD')=\op{Id}_{T(\HD)}$.
\endroster
\enddefinition

\subhead {3.1.3} \endsubhead
Let $i:Y\hookrightarrow \hx$ be the inclusion of a locally closed union
of simplices into $\hx$. 
By Lemma 2.2.7 we have $\psp\circ i^*=i^*\circ\psp$. 
Hence, the following diagram commutes 
(here $i^*:\Pre\L(\hx)\to\Pre\L(Y)$ is 
defined in the obvious way):
$$
\CD 
\SH_c(\hx) @<{\psp}<< \Pre\L(\hx) \\
@V{i^*}VV         @V {i^*}VV \\
\SH_c(Y) @<{\psp}<< \Pre\L(Y). \\ 
\endCD\tag{$*$}
$$
Notice that $Y$ is an open union of simplices in the
simplicial complex $\ov Y$, hence the functor
$\psp:\Pre\L(Y)\to\SH_c(Y)$ is obtained in the obvious way from
the functor $\psp:\Pre\L(\ov Y)\to\SH_c(\ov Y)$. 

Notice that each perverse simplex is a locally closed
union of simplices of $\hx$. Thus for each $\T\in\SS(X,\d)$
we have $\psp(\T)\in\SH_c(X,\d)$.

The restriction of $\psp$ to $\SS(X,\d)$ will be 
denoted by $\upp:\SS(X,\d)\to\SH_c(X,\d)$.  

The commutativity of the diagram $(*)$ is equivalent 
to the commutativity of the following diagram:
$$
\CD 
\SH_c(\hx) @>{\psm}>> \Pre\L(\hx) \\
@V{i^*}VV         @V {i^*}VV \\
\SH_c(Y) @>{\psm}>> \Pre\L(Y), \\ 
\endCD\tag{$**$}
$$
which means that there exists an isomorphism of functors
$\psm\circ i^*\overset\sim\to\to i^*\circ\psm$.

\subhead {3.1.4} \endsubhead
Let $\S\in\SH_c(X,\d)$. We will construct a functor
$\upm:\SH_c(X,\d)\to\SS(X,\d)$ and an isomorphism
$f:\upm\to\psm$ (here $\psm$ is restricted to $\SH_c(X,\d)$).
Let $\upm(\S)=\T$, $\psm(\S)=\T'$ and let
$Y=\PD$ be a perverse simplex. We set
$$
\T|_Y=\psm(\S|_{Y}).
$$
(Here $\psm(\S|_{Y}):=\psm(\S|_{\ov Y})|_Y$.)
The isomorphism of functors 
$\psm\circ i^*\overset\sim\to\to i^*\circ\psm$
provided by the diagram 3.1.3({$**$})
gives us the isomorphism for all $Y$
$$
f|_Y:\T|_Y\overset\sim\to\to\T'|_Y.\tag{$*$} 
$$
This is equivalent to providing isomorphisms 
$f(\D):T(\D)\to T'(\D)$ for all $\D\subeq\hx$,
commuting with restriction maps for all $Y$.
Let us complete the definition of $\T$. Let 
$\D_1\subeq\D_2$. We define
$$
t(\D_1,\D_2)=f^{-1}(\D_2)t'(\D_1,\D_2)f(\D_1).
$$
Notice that if $\D_1,\D_2\subeq Y$ for some perverse simplex $Y$, 
then $t(\D_1,\D_2)=\Id$ due to $(*)$. The definition of $\upm$
on morphisms is left to the reader. The stalkwise
isomorphisms $f(\D): T(\D)\to T'(\D)$ commuting with
restriction maps provide the isomorphism $\T\to\T'$.
Therefore the functors $\upm\overset\sim\to\to\psm$
are isomorphic.

\proclaim{3.1.5. Lemma} 
\roster
\item $\upm\circ\upp\simeq\Id$.
\item $\upp\circ\upm\simeq\Id$.
\endroster
\endproclaim

\demo{Proof} 
\roster
\item Clear, cf. 2.1.4, 2.2.8.
\item $\Id\simeq\psp\circ\psm\simeq\psp\circ\upm
\simeq\upp\circ\upm$, the first equivalence is Lemma 2.2.9.
\endroster
\enddemo

\subhead 3.2. Proof of the Theorem B2\endsubhead

\subhead {3.2.1}\endsubhead
Let $F:\Pre\L(\hx)\to\RR(\hx)$ be the isomorphism 
of 1.2. By abuse of notation we denote the full
subcategory of $\RR(\hx)$, corresponding to 
$\SS(X,\d)\sub\Pre\L(\hx)$ under $F$, by $\SS(X,\d)$.
In other words,

\definition{Definition} 
The category $\SS (X,\d)$
is a full subcategory of the abelian category $\RR (\hx)$. 
An object $\T$ of $\RR (\hx)$
belongs to $\SS(X,\d)$ 
if for any $\HD,\HD'\subseteq {^{-\d}\!\D}$ we have:
\roster
\item $T(\HD)=T(\HD')$,
\item if $\HD$ is a $\op{codim} 1$ face of $\HD'$, then
$t(\HD,\HD')=\op{Id}_{T(\HD)}$.
\endroster
\enddefinition

Until the end of this section $\geq$ stands for the relation in 
$\L(X,\d)$.

\proclaim{3.2.2. Lemma} 
Let $\hat\D,\hat\D'$ of $\hat X$
be such that  $\hat\D$ is a $\op{codim}\ 1$ face of $\hat\D'$.
Let $\HD\subseteq{^{-\d}\!\D}$,
$\hat\D'\subseteq{^{-\d}\!\D'}$. Then
$\D\lrar\D'$ and $\D\geq\D'$.
\endproclaim
\demo{Proof} Let
$\hat\D=\{v_1,\dots,v_k\}$, $\hat\D'=\{v_1,\dots,v_k,e\}$.
We will consider two cases.

Case 1. $e=\max(\HD')$.  Then $e$ is the barycenter of $\D'$.
By definitions $-\d(\D')>-\d(\D)$ i.e. $\d(\D)>\d(\D')$.
Let $v_i=\max(\HD)$. Then $v_i$ is the barycenter of $\D$. The
1-simplex $\{v_i, e\}$ is a simplex in $\hx$, thus $\D\lrar\D'$.
Moreover, $\D>\D'$ by Lemma 1.2.3.

Case 2. $e\neq \max(\hat\D')$. Then $v_i=\max(\HD)=\max(\HD')$,
and $v_i$ is the barycenter of both $\D$ and $\D'$, thus $\D=\D'$.
\enddemo

\proclaim{3.2.3. Theorem} The following categories
$$
\RR (X,\d)=\SS (X,\d)
$$
are isomorphic.
\endproclaim
\demo{Proof} The functor $\Phi:\Pre\L(X,\d)\to \SS(X,\d)$.
If $\S$ is an object of $\Pre\L(X,\d)$, then
$\T=\Phi(\S)$ is constructed as follows:
$$
T(\HD)=S(\D) \text {   for } \HD\subseteq\PD. 
$$
Let $\HD'$ be a $\op{codim}\ 1$ face of $\HD''$ 
and let $\HD'\subseteq\PD'$ and $\HD''\subseteq\PD''$.
By Lemma 3.2.2 $\D'\geq\D''$. We set:
$$
t(\HD', \HD'')=s(\D', \D'').
$$
Now let $\HD'\subeq\PD'$, $\HD_1\subeq\PD_1$,
$\HD_2\subeq\PD_2$, $\HD''\subeq\PD''$ be such
a quadruple for which we have to check the equivalence axiom. 
Then by Lemma 3.2.2 we have:
$$
\matrix
\D'\geq\D_1\geq\D''\ , \\
\D'\geq\D_2\geq\D''\ . \\
\endmatrix
$$
By definitions  we have:
$$
\aligned
t(\HD_1, \HD'')\circ t(\HD', \HD_1)
= &\ s(\D_1,\D'')\circ s(\D',\D_1) \\ 
= &\ s(\D',\D'')\\
= &\ s(\D_2,\D'')\circ s(\D',\D_2) \\
= &\ t(\HD_2, \HD'')\circ t(\HD', \HD_2).
\endaligned
$$
Let $h:\S_1\to\S_2$ be a morphism in $\Pre\L(X,\d)$.
We set:
$$
\Phi(h)(\HD)=h(\D) \text {   for } \HD\subseteq\PD. 
$$
Let $\T_1=\Phi(\S_1)$ and $\T_2=\Phi(\S_2)$.
Let $\HD'$ is a $\op{codim}\ 1$ face of $\HD''$
and let $\HD'\subseteq\PD'$ and $\HD''\subseteq\PD''$.
The commutativity of the diagram
$$
\CD
\T_1(\HD') @>{\Phi(h)(\HD')}>>\T_2(\HD') \\
@V{t_1(\HD',\HD'')}VV @V{t_2(\HD',\HD'')}VV \\
\T_1(\HD'') @>{\Phi(h)(\HD'')}>>\T_2(\HD'') \\  
\endCD
$$
follows from the commutativity of the diagram
$$
\CD
\S_1(\D') @>{h(\D')}>>\S_2(\D') \\
@V{s_1(\D',\D'')}VV @V{s_2(\D',\D'')}VV \\
\S_1(\D'') @>{h(\D'')}>>\S_2(\D''). \\  
\endCD
$$
\medskip

The functor $\Psi : \SS (X,\d)@>>>\RR (X,\d)$.
If $\T$ is an object of $\SS (X,\d)$, then
$\S=\Psi(\T)$ is constructed as follows:
$$
S(\D)=T(\HD) \text {   for } \HD\subseteq\PD. 
$$
$S(\D)$ is well defined since $\T$ is constant along 
$\PD$. 
Let $\D'$ and $\D''$ be two incident simplices of $X$
such that $\d(\D')=\d(\D'')+1$. We have to construct
the map $s(\D',\D'')$. Let $c'$ be the barycenter of $\D'$
and $c''$ be the barycenter of $\D''$. We set:
$$
s(\D',\D'')=t(\{c'\},\{c',c''\}),
$$
where $\{c',c''\}$ is a simplex in $\hx$.

We have to check the equivalence axiom.
Let $\D'>\D_b>\D''$, $\d(\D')=\d(\D_b)+1=\d(\D'')+2$.
We will consider two cases.

Case 1. $\D'$ and $\D''$ are incident. 
Let $c'$, $c_b$ and $c''$ be the barycenters of $\D'$, $\D_b$
and $\D''$ respectively. For the purposes of this proof we set
$s(\D',\D'')=t(\{c'\},\{c',c''\})$. 
Let us assign special names to
the following four simplices of $\hx$: 
$l=\{c',c''\}$, $l_b'=\{c',c_b\}$, $l_b''=\{c'',c_b\}$,
and $tr=\{c',c_b, c''\}$. We have:
$$
\aligned
s(\D',\D'')= &\ t(c',l)\\ 
= &\ t(l, tr)\circ t(c',l) \\ 
= &\ t(l_b', tr)\circ t(c', l_b') \\
= &\ t(l_b', tr)\circ t(c_b,l_b')\circ t(c',l_b') \\  
= &\ t(l_b'', tr)\circ t(c_b,l_b'')\circ t(c',l_b')\\ 
= &\ t(c_b,l_b'')\circ t(c',l_b')\\
= &\ s(\D_b,\D'')\circ s(\D',\D_b). \\
\endaligned\tag {$*$}
$$
The equivalence axiom follows.

Case 2. $\D'$ and $\D''$ are not incident. 
By Lemma 1.2.3 $\d(\D')=1$, $\d(\D'')=-1$. Then
there is exactly one $\D_0$ such that 
$\D'\lrar\D_0\lrar\D''$, $\d(\D_0)=0$, and the
equivalence axiom is vacuous.

Let $h:\T_1\to\T_2$ be a morphism in $\SS(X,\d)$.
We set:
$$
\Psi(h)(\D)=h(\HD) \text {   for } \HD\subseteq\PD. 
$$
$\Psi(h)(\D)$ is well defined since $\T_1$ and $\T_2$
are constant along $\PD$.
Let $\S_1=\Psi(\T_1)$ and $\S_2=\Psi(\T_2)$, and
let $\d(\D)=\d(\D')+1$, and $\D\lrar\D'$. 
The commutativity of the diagram
$$
\CD
\S_1(\D) @>{\Psi(h)(\D)}>>\S_2(\D) \\
@V{s_1(\D,\D')}VV @V{s_2(\D,\D')}VV \\
\S_1(\D') @>{\Psi(h)(\D')}>>\S_2(\D') \\  
\endCD
$$
follows from the commutativity of the diagram
$$
\CD
\T_1(\{c\}) @>{h(\{c\})}>>\T_2(\{c\}) \\
@V{t_1(\{c\},\{c,c'\})}VV @V{t_2(\{c\},\{c,c'\})}VV \\
\T_1(\{c,c'\}) @>{h(\{c,c'\})}>>\T_2(\{c,c'\}). \\  
\endCD
$$
\medskip

Let $G: \RR(X,\d)\to\Pre\L(X,\d)$ be the isomorphism 
of 1.2. It follows from our explicit construction that
$\Psi\circ\Phi\circ G=\Id$.
\medskip 

$\Phi\circ G\circ\Psi=\Id$. Let $\T\in\SS(X,\d)$,
and let $\Phi\circ G\circ\Psi(\T)=\check\bold T$.
For a simplex $\HD$ in $\hx$:
$$
\check T(\HD)=T(\HD).
$$
Let $\HD'$ be a $\op{codim}\ 1$ face of $\HD''$ 
and let $\HD'\subseteq\PD'$ and $\HD''\subseteq\PD''$. 
By Lemma 3.2.2 $\D'\lrar\D''$ and $\D'\geq\D''$.
If $\D'=\D''$, then 
$\check t(\HD',\HD'')=\Id=t(\HD',\HD'')$.
If $\D'>\D''$, then there exists a sequence
$\D'=\D_0, \D_1, \D_2,\dots,\D_m=\D''$ such that:
\roster
\item $\D_i\lrar\D_{i+1}\lrar\D_{i+2}$, $0\leq i\leq m-2$,
\item $\d(\D_i)=\d(\D_{i+1})+1$, $0\leq i\leq m-1$,
\item $\D'\lrar\D_i\lrar\D''$, $0\leq i\leq m$ (Lemma 1.2.3).
\endroster
Let $c_i$ be the barycenter of $\D_i$. If 
$\S=G\circ\Psi(\T)$, then by definitions 
$s(\D_i,\D_{i+1})=t(\{c_i\},\{c_i,c_{i+1}\})$.
Then
$$
\aligned
\check t(\HD',\HD'')=&\ s(\D',\D'')\\
=&\ s(\D_{m-1},\D_m)\circ\dots\circ s(\D_0,\D_1)\\
=&\ t(\{c_{m-1}\},\{c_{m-1},c_m\})\circ\dots\circ
t(\{c_0\},\{c_{0},c_1\})
\endaligned
$$
We claim that: 
$$
t(\{c_0\},\{c_{0},c_m\})=
t(\{c_{m-1}\},\{c_{m-1},c_m\})\circ\dots\circ
t(\{c_0\},\{c_{0},c_1\})=\check t(\HD',\HD'').
$$
If $m=1$ then there is nothing to prove. In general the claim
is proved by induction. The step of induction is an
argument similar to ($*$).
Since $\HD'\subseteq\PD'$ and $\HD''\subseteq\PD''$,
$c'=\max\HD'$, and $c''=\max\HD''$ 
($c'$ and $c''$ are barycenters of $\D'$ and $\D''$ respectively).
We have:
$$
\matrix
\{c'\}\subeq\HD'\sub\HD'' \\
\{c'\}\sub \{c',c''\}\subeq\HD''
\endmatrix
$$
Notice that $c'\subeq\PD'$, $\{c',c''\}\subeq\PD''$.
Let $\T'\in\SS(X,\d)\sub\Pre\L(\hx)$ be the 
presheaf corresponding to $\T\in\SS(X,\d)\sub\RR(\hx)$.
We have:
$$
\aligned
t(\{c'\},\{c',c''\}) 
= &\ t'(\{c'\},\{c',c''\}) \\
= &\ t'(\{c',c''\},\HD'')\circ t'(\{c'\},\{c',c''\})\\
= &\ t'(\{c'\},\HD'') \\
= &\ t'(\HD',\HD'')\circ t'(\{c'\},\HD') \\
= &\ t'(\HD',\HD'') \\
= &\ t(\HD',\HD'')
\endaligned
$$
Therefore $\check t(\HD',\HD'')=t(\HD',\HD'')$ and 
$\check\bold T=\T$.
\enddemo

\subhead {3.2.4} \endsubhead
Summarizing, we have constructed the equivalence functors:
$$
\SH_c(X,\d)@>{\upm}>>\SS(X,\d)=\SS(X,\d)@>{\Psi}>>\RR(X,\d).
$$
Recall that in 1.2 we constructed the isomorphism functor:
$$
\Xi_B:\RR(X,\d)\to\modb(X,\d).
$$
Let $\Th(X,\d):\SH_c(X,\d)\to\modb(X,\d)$ be the 
equivalence functor defined as 
$\Th(X,\d)=\Xi_B\circ\Psi\circ\upm$.
In other words, we have

\proclaim{Theorem B2} The functor
$$
\Th(X,\d):\SH_c(X,\d)@>{\sim}>>\modb(X,\d)   
$$
is an equivalence of categories.
\endproclaim

\subhead 3.3. Injective objects in $\RR(X,-\d)$
\endsubhead

In this section we list some simple facts concerning injective
objects in the category $\RR(X,-\d)$. The material 
of this section is a generalization of some of the facts 
presented in section 2.4.

\subhead{3.3.1}\endsubhead
In this subsection we construct injective objects in 
the category $\Pre\L(X,-\d)$. Our construction generalizes 
injective simplicial sheaves (see 2.4). Let $\leq$
denote the relation in $\L(X,\d)$.

\definition{Definition} Let $\D\sub X$ be a simplex
and let $V$ be a vector space. We define the object
$[{}^{\d}\!\D]^V$ as follows:
$$
[{}^{\d}\!\D]^V(\D')=\cases
V, & \D'\leq\D \\
0, & \D'\not\leq\D, \\
\endcases
$$
where the restriction maps between copies of $V$ are all
identity maps. We will also write 
$[{}^{\d}\!\D]=[{}^{\d}\!\D]^\F$.
It is easy to see that for any $\D$ and $V$,
 $[{}^{\d}\!\D]^V$ is an injective object in $\Pre\L(X,-\d)$. 
\enddefinition

\proclaim{3.3.2. Lemma} Let $\S\in\Pre\L(X,-\d)$, and let
$[{}^{\d}\!\D]^V$ be an injective object. Then
\roster
\item $\Hom(\S,[{}^{\d}\!\D]^V)=\Hom_\F (S(\D),V)$.
\item In particular, if $\D''\leq\D'$, then
$\Hom([{}^{\d}\!\D']^V,[{}^{\d}\!\D'']^W)=\Hom_\F(V,W)$.
\endroster
\endproclaim

\demo{Proof} The proof is left to the reader.
\enddemo
 
\subhead{3.3.3}\endsubhead
In 1.2 we have seen that $\Pre\L(X,-\d)=\RR(X,-\d)$. 
We transport the injective objects $[{}^{\d}\!\D]^V$  
to the category $\RR(X,-\d)$ preserving the
notation.

Let $\FF$ be the ``constant sheaf'' in $\RR(X,-\d)$, i.e.:
$$
F(\D)=\F \qquad \text{ for all } \D\in X,
$$ 
where the restriction maps between copies of $\F$
are all identity maps.

\definition{Definition} Let $\S\in\RR(X,-\d)$. The
vector space $\G(X;\S)$ of {\sl global sections} of $\S$
is defined as follows: 
$$
\G(X;\S)=\Hom_{\RR(X,-\d)}(\FF,\S).
$$
\enddefinition

\proclaim{Lemma} Let $\D\sub X$ be a simplex
and let $V$ be a vector space. Then
$$
\G(X;[{}^{\d}\!\D]^V)=V.
$$ 
\endproclaim

\demo{Proof} Indeed, by definitions and Lemma 3.3.2:
$$
\G(X;[{}^{\d}\!\D]^V):=\Hom_{\RR(X,-\d)}(\FF,[{}^{\d}\!\D]^V)
=\Hom_{\F}(\F,V)=V.
$$
\enddemo

\head
{Chapter 4. 
Koszul duality and perverse sheaves}
\endhead

\subhead 4.1. Mixed quiver algebras \endsubhead
\definition{4.1.1. Definition} (cf. [BGSo, 4.1]) A quiver $\Q$
is called {\sl mixed} if it is equipped with a function
$w:V(\Q)\to\Bbb Z$ (called a {\sl weight}) such that for any
two vertices $\a,\b\in V(\Q)$ there is no arrow 
$(\a\to\b)\in E(\Q)$
if $w(\a)\leq w(\b)$.
\enddefinition

Notice that this definition is different from the 
one given in [Vyb1], where we also assumed that: 

$(*)$ For $\a,\b\in V(\Q)$ there is no arrow $(\a\to\b)\in E(\Q)$ 
unless $w(\a)=w(\b)+1$.\newline
In the rest of this section we will consider only mixed quivers
satisfying $(*)$.  

Note that the algebra $\F\Q$ has a natural 
grading:
$$
\F\Q_m=\{\text {vector space spanned by all paths of length } m\}
$$
We denote by $\F\Q_+=\oplus_{m>0}\F\Q_m$. 

\definition{4.1.2. Definition} {\sl A mixed quiver algebra}
$C=\F\Q/J$ is the quotient of $\F\Q$ by a homogeneuos ideal
$J\subset\F\Q_+$.
\enddefinition

A mixed quiver algebra $C$ inherits its grading from $\F\Q$ since 
$J$ is a homogeneous ideal. Let $I$ be the set of local idempotents 
$e(\a)$ with $\a\in V(\Q)$. Since there is a one-to one correspondence
$I=V(\Q)$ we may assume that $w$ is a function on $I$,
$w:I\to\Bbb Z$. We define $I_l:=w^{-1}(-l)$.

\subhead {4.1.3}\endsubhead
In this subsection we consider the category $\modc$ of left finitely 
generated modules over a mixed quiver algebra $C$. 
If $M$ is such a module,
$M$ admits a natural (standard) grading:                                            
$$
M_l=\sum_{e\in I_l}eM, \qquad
\text{and} \qquad 
M=\oplus_{l} M_l.
$$
We can shift 
the standard grading by any integer $k$:
$(M\langle k\rangle)_l=M_{l-k}$ ($k$-grading).
In this terminology standard grading is $0$-grading.

Now let us consider the category $\modgrc$ of finitely generated
graded modules over the graded algebra $C$. If $M$ is an object of
$\modgrc$ and $x\in M_i$ is a homogeneous element of $M$,
then $\op{grdeg}(x)=i$ denotes the homogeneous degree of $x$.
Note that both $\modc$ and $\modgrc$ are abelian categories
since $C$ is Noetherian.

\proclaim{4.1.4. Proposition [Vyb1]} The abelian category $\modgrc$ 
is the direct sum:
$$
\modgrc=\bigoplus_{k\in\Bbb Z} (\modc)_k ,
$$
where $(\modc)_k$ is (isomorphic to) the category $\modc$ with  
$k$-grading.
\endproclaim

\proclaim{Corollary} The bounded derived category
$\DD(\modgrc)$ is the direct sum:
$$
\DD(\modgrc)=\bigoplus_{k\in\Bbb Z}\DD ((\modc)_k).
$$
\endproclaim

\remark{Remark} A similar situation
is discussed in [PP].
\endremark

\subhead{4.1.5}\endsubhead
Let $X$ be a connected finite simplicial complex, and
let $\d$ be a perversity. Notice that the quiver $\Q(X,\d)$
is mixed with weigt function $w=\d$. Hence the algebras
$A(X,\d)$ and $B(X,\d)$ introduced in 1.2 
are mixed quiver algebras. The algebras $A(X,\d)$ and $B(X,\d)$
are also quasi-hereditary and of finite global dimension [CPS, Vyb1].

\subhead 4.2. Koszul duality \endsubhead
In this section we complete the proof of Theorem C,
Theorem D1, and Theorem D2. The material is mostly
borrowed from [BGSc, BGSo, Mac2, Mac3, Pol, Vyb 1, Vyb2].

\definition{4.2.1. Definition [BGSo]}
A positively graded ring $C=\oplus_{j\geq 0}C_j$ 
with semisimple $C_0$ is called Koszul if $C_0$ considered
as a graded left $C$-module admits a graded projective resolution
$$
\dots@>>> P^2 @>>> P^1 @>>> P^0 @>>> C_0 @>>> 0
$$
such that $P^i$ is generated by its component of degree $i$,
$P^i=CP^i_i$.
\enddefinition

\proclaim{Proposition [BGSo, 2.1.3]} Let $C=\oplus_{j\geq 0}C_j$ 
be a positively graded ring and suppose $C_0$ is semisimple.
The following conditions are equivalent:
\roster
\item $C$ is Koszul.
\item For any two pure $C$-modules $M$, $N$ of weights $m$, $n$
respectively we have $\op{Ext}_{\modgrc}^i(M,N)=0$
unless $i=m-n$.
\item $\op{Ext}_{\modgrc}^i(C_0,C_0\langle n\rangle)=0$
unless $i=n$.
\endroster
\endproclaim

\subhead{4.2.2}\endsubhead
Until the end of this section $X$ is a finite connected simplicial 
complex, and $\d$ is a perversity.

\proclaim{Lemma}The following conditions are equivalent:
\roster
\item $B(X,\d)$ is Koszul.
\item For any two simple $B(X,\d)$-modules $S_e$ and $S_{e'}$
we have $\op{Ext}_{B(X,\d)}^i(S_e,S_{e'})=0$
unless $i=\d(e)-\d(e')$.
\endroster
\endproclaim

There is an obvious analogous statement for $A(X,\d)$.  

\demo{Proof} By Corollary 4.1.4 we have 
$$
\op{Ext}^i_{\modgrb(X,\d)}(S_e,S_{e'})=
\op{Ext}^i_{B(X,\d)}(S_e,S_{e'}),
$$
where we consider
$S_e$ and $S_{e'}$ as graded modules equipped with standard 
grading. Then it is easy to see that the condition (2)
is equivalent to the condition (2) of Proposition 4.2.1.
\enddemo

\proclaim{4.2.3. Lemma} 
\roster
\item The (left) quadratic
dual of the quadratic algebra $A(X,\d)$ 
is isomorphic to $B(X,-\d)$, $A(X,\d)^!=B(X,-\d)$.
\item The (left) quadratic
dual of the quadratic algebra $B(X,\d)$ 
is isomorphic to $A(X,-\d)$, $B(X,\d)^!=A(X,-\d)$.
\endroster
\endproclaim

\demo{Proof} The proof is straightforward since
the chain complex relations,
and the equivalence relations 
are quadratic dual to each other.
\enddemo

\subhead{4.2.4}\endsubhead
Let $\M(X,\d)$ be the category of constructible perverse sheaves,
and let $\S_{\D}\in\M(X,\d)$ be the simple perverse sheaf
associated to a simplex $\D$ (see 1.3).
Lemma A and Lemma C below are borrowed from [Pol].

\proclaim{Lemma A [Pol, 1.2]} Let $\D$ and $\D'$ be 
two simplices in $X$. Then
$$
\op{Ext}^i_{\DD_c(X)}(\S_{\D},\S_{\D'})=
\cases
\F, & \D\geq\D' \text{ and } i=\d(\D)-\d(\D') \\
0, & \text{otherwise}. \\
\endcases
$$
\endproclaim

\proclaim{Lemma B} The algebra $B(X)$ is Koszul.
\endproclaim

\demo{Proof} The following categories 
$\DD_c(X)\simeq\DD(\SH_c(X))\simeq\DD(\modb(X))$ are 
equivalent (see 1.3 and 2.2). Let $e$ and $e'$ be the local 
idempotents of $B(X)$ corresponding to simplices $\D$ and $\D'$.
Then
$$
\op{Ext}^i_{B(X)}(S_e,S_{e'})=
\op{Ext}^i_{\DD(\SH_c(X))}(\S_\D,\S_{\D'})=
\op{Ext}^i_{\DD_c(X)}(\S_\D,\S_{\D'})=0
$$ 
unless
$i=\d(\D)-\d(\D')=\dim\D'-\dim\D=\d(e)-\d(e')$.
(In this case $\d=\text{bottom perversity}$.)
Then $B(X)$ is Koszul by Lemma 4.2.2.
\enddemo

\proclaim{Lemma C [Pol, 2]} The following two algebras
$$
\oplus_{\D,\D'}\op{Ext}^*_{\DD_c(X)}(\S_{\D},\S_{\D'})=B(X,\d)
$$
are canonically isomorphic.
\endproclaim

\subhead{4.2.5}\endsubhead
In this subsection we sketch some material mostly borrowed 
from [BGSo, Pol, Vyb2].
Theorem C below (in a slightly different form) is due to R. 
MacPherson [Mac2, Mac3].

Any object $\S\in\M(X,\d)$ has a canonical finite
increasing filtration $W_\bu=W_\bu\S$ such that
$\op{gr}^W_i\S=W_i\S/W_{i-1}\S$ is a direct sum of simple
objects $\S_\D$ with $\d(\D)=i$ 
(see [BBD, 5.3.6], [BGSo, 4.1.2], [Pol, 4.1]).
Let us define a functor: 
$$
\matrix 
\Sigma^i(X,\d):\M(X,\d)\to\Vect,\\
\Sigma^i(X,\d)(\S)=
\Hom_{\M(X,\d)}(\oplus_{\D}\S_\D,\op{gr}^W_{-i}\S).
\endmatrix
$$
The functor 
$$
\matrix
\Sigma(X,\d)=\oplus_i\Sigma^i(X,\d):\M(X,\d)\to\op{Vectgr}_\F,\\
\Sigma(X,\d)(\S)_i=\Sigma^i(X,\d)(\S)
\endmatrix
$$
assignes a graded vector space to each $\S\in\M(X,\d)$.
Let us also define:
$$
\ti A_i=\oplus_a\Hom(\Sigma^a(X,\d),\Sigma^{a+i}(X,\d)).
$$
The composition of morphisms defines a product 
$\ti A_i\times\ti A_j\to\ti A_{i+j}$. Note that
$\ti A_0=\op{End}(\oplus_{\D}\S_\D)$. The graded
algebra $\ti A$ acts canonically on the graded vector space
$\Sigma(X,\d)(\S)$ for any $\S\in\M(X,\d)$. One can show that the
functor $\Sigma(X,\d):\M(X,\d)\to\op{mod-}\!\ti A$
is an equivalence of categories, and  
that
the algebras, $\ti A=B(X,-\d)^!=A(X,\d)$ are isomorphic 
(see [Pol, 4]). 
In other words, we have
 
\proclaim{Theorem C} The functor
$$
\Sigma(X,\d):\M(X,\d)@>{\sim}>>\moda(X,\d)
$$
is an equivalence of categories.
\endproclaim

\subhead{4.2.6}\endsubhead
The following theorem is borrowed from [Pol].
The statement (3) was obtained by R. MacPherson [Mac3].

\proclaim{Theorem} 
\roster
\item The algebra $B(X,\d)$ is Koszul.
\item The algebra $A(X,\d)$ is Koszul.
\item The following triangulated categories 
$$
\DD(\M(X,\d))\simeq\DD_c(X)
$$
are equivalent.
\endroster
\endproclaim

\remark{Remark} A. Polishchuk shows in [Pol] that (3) follows from
(2). Using Lemma 4.2.2, Lemma 4.2.4.A, and Theorem C, 
it is also possible to show that (2) follows from (3).
\endremark

\subhead{4.2.7}\endsubhead
Here we adapt the construction of the Koszul duality functors
from [BGSo, 2.12] to our situation. 
We denote $k=A(X,\d)_0=B(X,\d)_0$.

The functor 
$K:\CC(\modgra(X,\d))\to\CC(\modgrb(X,-\d))$
is constructed as follows. Let $M=M^\bu\in\CC(\modgra(X,\d))$.
Then
$$
\matrix
KM=B(X,-\d)\otimes_k M, \\
(KM)^p_q=\bigoplus_{p=i+j,q=l-j}B(X,-\d)_l\otimes_k M^i_j .
\endmatrix
$$
Similarly, the functor 
$K':\CC(\modgra(X,\d))\to\CC(\modgrb(X,-\d))$
is constructed as follows. Let $N=N^\bu\in\CC(\modgra(X,\d))$.
Then
$$
\matrix
K'N=\Hom_k (B(X,-\d), M), \\
(K'N)^p_q=\bigoplus_{p=i+j,q=l-j}\Hom_k(B(X,-\d)_{-l},M^i_j) .
\endmatrix
$$
The differentials of $KM^\bu$ and $K'N^\bu$ are defined in 
[BGSo, 2.12.1].
Since the functors $K$ and $K'$ are exact they induce the
functors 
$\DD(K), \DD(K'):\DD(\modgra(X,\d))\to\DD(\modgrb(X,-\d))$.

\proclaim{Theorem [BGSo, 2.12.6]} 
\roster
\item The functor
$$
\DD(K):\DD(\modgra(X,\d))@>{\sim}>>\DD(\modgrb(X,-\d))
$$
is an equivalence of triangulated categories.
\item The functor
$$
\DD(K'):\DD(\modgra(X,\d))@>{\sim}>>\DD(\modgrb(X,-\d))
$$
is an equivalence of triangulated categories.
\endroster
\endproclaim

Recall that by Corollary 4.1.4: 
$$
\split
\DD(\modgra(X,\d))=\bigoplus_{m\in\Bbb Z}\DD((\moda(X,\d))_m),\\
\DD(\modgrb(X,\d))=\bigoplus_{m\in\Bbb Z}\DD((\modb(X,\d))_m).\\
\endsplit
$$
\proclaim{Lemma} 
\roster
\item Let $M\in\DD((\moda(X,\d))_0)$. Then
$\DD(K)(M)\in\DD((\modb(X,-\d))_0)$.
\item Let $N\in\DD((\moda(X,\d))_0)$. Then
$\DD(K')(N)\in\DD((\modb(X,-\d))_0)$.
\endroster
\endproclaim

\demo{Proof} The proof is left to the reader.
\enddemo

\proclaim{Corollary}The functors
$$
\DD(K),\ \DD(K'):\DD(\moda(X,\d))@>{\sim}>>\DD(\modb(X,-\d))
$$
are equivalences of triangulated categories.
\endproclaim

\subhead{4.2.8}\endsubhead
The category $\DD(\moda(X,\d))$ has the standard 
$t$-structure with the core $\moda(X,\d)$. The
Koszul duality functors $K$ and $K'$ transform
this standard $t$-structure to non-standard
$t$-structures in $\DD(\modb(X,-\d))$ (see [BGSo, 2.13]).
We will give a description of the cores of these
non-standard $t$-structures as follows.

\proclaim{Lemma} (cf. [BGSo, 2.13.3])
\roster
\item Let $M\in\moda(X,\d)$. Then 
$KM\in\DD(\modb(X,-\d))$ is a complex of
projective $B(X,-\d)$-modules
$$
\dots @>>>P^i@>>>P^{i+1}@>>>\dots
$$
such that for any $i$, $P^i$ is generated by its component
of degree $-i$, $P^i=B(X,-\d)P^i_{-i}$.
\item Let $N\in\moda(X,\d)$. Then 
$K'N\in\DD(\modb(X,-\d))$ is a complex of
injective $B(X,-\d)$-modules
$$
\dots @>>>I^i@>>>I^{i+1}@>>>\dots
$$
such that for any $i$, $I^i$ is cogenerated by its component
of degree $-i$, $I^i=\Hom_k(B(X,-\d),I^i_{-i})$. 
\endroster
\endproclaim

\demo{Proof} 
\roster
\item By construction we have $P^j=B(X,-\d)\otimes_k M_j$, which
is a projective module generated by $M_j$.
($M_j$ is considered as a semisimple $B(X,-\d)$-module
concentrated in degree $-j$.)
\item By construction we have $I^j=\Hom_k(B(X,-\d),N_j)$, which
is an injective module cogenerated by $N_j$.
\endroster
\enddemo

\subhead{4.2.9}\endsubhead
Since the categories $\moda(X,\d)=\PP(X,\d)$
(resp. $\modb(X,\d)=\RR(X,\d)$) 
are isomorphic, we can transport the Koszul duality
functor $K':\CC(\moda(X,\d))@>>>\CC(\modb(X,-\d))$
to the functor 
$\bold K':\CC(\PP(X,\d))@>>>\CC(\RR(X,-\d))$ (cf. 4.3).
Let $\S\in\PP(X,\d)$, and let us consider $\PP(X,\d)$ as the
core of the standard $t$-structure on $\DD(\PP(X,\d))$.
We will describe $\bold K'\S\in\CC(\RR(X,-\d))$ explicitly
as follows: 
$$
\bold {K'S}^i=\bigoplus_{\d(\D)=-i}[{}^{\d}\!\D]^{S(\D)}\ .
$$
The differential
$$
\matrix\format\l \\
\bold {K'S}^i@>{d^{i}}>>\bold {K'S}^{i+1}  \\
\bigoplus_{\d(\D)=-i}
[{}^{\d}\!\D]^{S(\D)}@>{d^{i}}>>
\bigoplus_{\d(\D)=-i-1}
[{}^{\d}\!\D]^{S(\D)}
\endmatrix
$$
is the zero map between components $[{}^{\d}\!\D]^{S(\D)}$ and 
$[{}^{\d}\!\D']^{S(\D')}$
if $\D'\not\le\D$, and is induced by 
$s(\D,\D'):S(\D)\to S(\D')$
if $\D'<\D$, and $\d(\D)=\d(\D')+1$.
(Here $\le$ refers to the relation in $\L(X,\d)$.)

Let us consider
the complex $\G(X;(\bold K'\S)^\bu)$ of global sections 
(see 3.3) of $(\bold {K'S})^\bu$:
$$
\G(X;\bold {K'S}^i)=\bigoplus_{\d(\D)=-i}{S(\D)}
$$
The differential 
$d^i:\G(X;\bold {K'S}^i)\to\G(X;\bold {K'S}^{i+1})$
is the zero map between $S(\D)$ and
$S(\D')$ if $\D'\not\le\D$, and is given by 
$s(\D,\D'):S(\D)\to S(\D')$
if $\D'<\D$, and $\d(\D)=\d(\D')+1$. We will denote
the cohomology spaces of this complex by $H^i(X;\bold K'\S)$.

\subhead{4.2.10}\endsubhead
Let $\S\in\PP(X,\d)$. The {\sl total chain complex} $C_\bu(X;\S)$
associated to $\S$ is constructed as follows [Mac2]: 
$$
C_i(X;\S)=\bigoplus_{\d(\D)=i}S(\D).
$$
The differential $d_i:C_i(X;\S)\to C_{i-1}(X;\S)$ is the
map whose matrix elements are $s(\D,\D'):S(\D)\to S(\D')$
if $\D'<\D$ and $\d(\D)=\d(\D')+1$, and zero otherwise.
The homology spaces
of $C_\bu(X;\S)$ are denoted by $H_i(X;\S)$.

\proclaim{Lemma} Let $\S\in\PP(X,\d)$. Then the complexes
$$
(C_i(X;\S), d_i)=(\G(X;\bold {K'S}^{-i}), d^{-i})
$$
coincide. Hence:
$$
H_i(X;\S)=H^{-i}(X;\bold K'\S).
$$
\endproclaim

\demo{Proof} This is clear by inspection.
\enddemo

\subhead{4.2.11}\endsubhead
We have constructed the equivalence functors:
$$
\CD
\DD(\M(X,\d)) &   & \DD(\SH(\TT(X,\d))) \\
@V{\DD(\Sigma(X,\d))}VV    @V{\DD(\Omega(X,-\d))}VV \\
\DD(\moda(X,\d))@>{\DD(K)}>>\DD(\modb(X,-\d)) \\
\endCD
$$
Let $\wti{K}_1:\DD(\M(X,\d))\to\DD(\SH(\TT(X,-\d)))$
be the equivalence functor defined as
$\wti{K}_1=
\DD(\Omega(X,-\d))^{-1}\circ\DD(K)\circ\DD(\Sigma(X,\d))$
In other words, we have

\proclaim{Theorem D1} The functor
$$
\wti{K}_1:\DD(\M(X,\d))@>{\sim}>>\DD(\SH(\TT(X,-\d))) 
$$
is an equivalence of triangulated categories.
\endproclaim

Similarly, we have:
$$
\CD
\DD(\M(X,\d)) &   & \DD(\SH_c(X,-\d)) \\
@V{\DD(\Sigma(X,\d))}VV    @V{\DD(\Th(X,-\d))}VV \\
\DD(\moda(X,\d))@>{\DD(K)}>>\DD(\modb(X,-\d)) \\
\endCD
$$
Let $\wti{K}_2:\DD(\M(X,\d))\to\DD(\SH_c(X,-\d))$
be the equivalence functor defined as
$\wti{K}_2=
\DD(\Th(X,-\d))^{-1}\circ\DD(K)\circ\DD(\Sigma(X,\d))$.
In other words, we have

\proclaim{Theorem D2} The functor
$$
\wti{K}_2:\DD(\M(X,\d))@>{\sim}>>\DD(\SH_c(X,-\d)) 
$$
is an equivalence of triangulated categories.
\endproclaim

Clearly, one could also state a version of Theorem D1
and Theorem D2 for the functor $K'$ simply by substituting 
$K'$ instead of $K$ everywhere.
The corresponding
equivalence functors are denoted by $\wti{K}'_1$ and
$\wti{K}'_2$ respectively. 
The latter functor is used in the following section.

\subhead{4.2.12. A corollary of the Theorem D2: perverse sheaves
as complexes of sheaves with allowable support}\endsubhead

Let us fix a perversity $\ov p$ as in [GM1], and the associated
cellular perversity $\d$.

\definition{Definition [GM1]}
\roster
\item A pseudomanifold of dimension $n$ is a compact space $X$ such that
$X$ is the closure of the union of the $n$-simplices in any triangulation
of $X$, and each $n-1$ simplex is a face of exactly two $n$-sipmlices.
\item A stratification of a pseodomanifold is a filtration by closed subspaces
$$
X=X_n\sup X_{n-1}=X_{n-2}\sup X_{n-3}\sup\dots\sup X_0
$$
satisfying certain conditions (see [GM1, 1.1]). We can represent $X$ 
as a disjoint union $X=\sqcup_{S\in\SS} S$ where the strata $S$ are 
connected components of
$X_i-X_{i-1}$.
\item If $i$ is an integer, a subspace $Y\sub X$ is called
{\sl $(\ov p, i)$-allowable} if $\dim (Y)\leq i$ and
$\dim (Y\cap X_{n-k})\leq i-k+\ov p(k)$ for all $k\geq 2$.
\endroster 
\enddefinition

Let $\TT$ be a triangulation of a pseudomanifold $X$ refining 
a stratification $\SS$.

Let $\Cal M_{\SS}(X,\d)$ be category of perverse 
sheaves constructible with respect
to the stratification $\SS$. Let $\M(X,\d)$ be our usual
category of perverse 
sheaves constructible with respect
to the triangulation $\TT$. 

The refinement functor (see 1.3.9)
provides the embedding 
$\refn: \Cal M_{\SS}(X,\d)\hookrightarrow \M(X,\d)$.
Let us take $\A\in\Cal M_{\SS}(X,\d)$ and consider $\wti K'_2(\refn\A)$
as the complex of sheaves constant along $\d$-perverse simplices.

\proclaim{Corollary D3} Let $\A\in\Cal M_{\SS}(X,\d)$.
For\ $i\in\Bbb Z$, $-\ov p(n)\leq i\leq n-\ov p(n)$, 
let $j=n-\ov p(n)-i$, and let $\ov q$ be 
another
perversity defined by $\ov q(k)=\ov p(n)-\ov p(n-k)$.
Then $\wti K'_2(\refn\A)$ is (quasi-isomorphic to) a complex of sheaves 
$$
\dots @>>> \B^{i-1} @>>> \B^{i} @>>> \B^{i+1}@>>> \dots
$$
on $X$ such that $\supp\B^i$ is $(\ov q,j)$-allowable with respect
to the stratification $\SS$.
\endproclaim

\demo{Proof} We will show that in fact the claim is  
true for any $\A\in\M(X,\d)$.
It is clear from our construction (in particular 4.2.8.(2)) that
$\supp\B^i\sub X^{\d}_{-i}$, and $\supp\B^i$ is a union of perverse simplices,
and therefore a union of sipmlices of the first barycentric subdivision of $\TT$.
By Lemma 1.1.9 we have
$Q^{\ov q}_j=X^{\d}_{\ov p(n)-n+j}$. Then 
$Q^{\ov q}_{j=n-\ov p(n)-i}=X^{\d}_{-i}$.
Since $Q^{\ov q}_{j}$ is $(\ov q,j)$ allowable with respect to the
stratification $\SS$ [GM1, 3.2], and 
$\supp\B^i\sub Q^{\ov q}_{j}$ we conclude that 
$\dim\supp\B^i\leq j$ and $\dim(\supp\B^i\cap X_{n-k})\leq j-k+\ov q(k)$,
i.e. $\supp\B^i$ is $(\ov q,j)$-allowable.
\enddemo

\remark{Remark} By the above Corollary the category of 
stratification-constructible 
perverse sheaves $\Cal M_{\SS}(X,\d)$ is embedded by the duality 
functor into the category
of complexes of sheaves with the allowable support. 
One may ask how to characterize $\SS$-constructible 
perverse sheaves inside this larger category. This very
interesting question will be addressed in [GMMV]. 
\endremark

\subhead 4.3. Verdier and Koszul duality for
simplicial sheaves and cosheaves
\endsubhead

In this section we study the relationship between the two 
dualities. $X$ is assumed to be a finite connected simplicial 
complex.

\subhead{4.3.1}\endsubhead
Recall that $B(X)=B(X,\text{bottom perversity})$,
and let us denote $A'(X)=A(X,\text{top perversity})$,
and $k=A'(X)_0=B(X)_0$. In 4.2 we have seen that 
$A'(X)$ is the (left) quadratic dual of $B(X)$,
$A'(X)=B(X)^!$. We have also seen that $A'(X)$
and $B(X)$ are Koszul algebras. We will consider
the quasi-inverse Koszul duality functors 
(which we denote by $L$ and $L^{-1}$ in this section)
for $A'(X)$ and $B(X)$, from a point of view  somewhat different 
from 4.2. The rest of this subsection is essentially
lifted from [BGSo, 2.12]. 

The functor $L:\CC(\modgrb(X))\to\CC(\modgra'(X))$
is defined as follows.
Let $M^\bu\in\modgrb(X)$. Then
$$
LM^\bu=A'(X)\otimes_k M^\bu.
$$
The differential of $LM^\bu$ is (up to signs) the sum of 
the differential coming from the complex $M^\bu$ and the Koszul complex, 
see [BGSo, 2.12.1].

The functor $L^{-1}:\CC(\modgra'(X))\to\CC(\modgrb(X))$
is defined as follows. 
Let $N^\bu\in\CC(\modgra'(X))$.
We define a bicomplex: 
$$
(L^{-1}N)^{i,j}=\Hom_k(B(X),N^i_{j}).\tag {$*$}
$$
The differentails of this bicomplex are constructed in 
[BGSo, 2.12.1] (also see above). The
diagonal complex of the bicomplex $(*)$ is $L^{-1}N^\bu$.

If $M^\bu$ is equipped with the standard grading, then
$LM^\bu$ is also equipped with the standard grading, and 
likewise, if $N^\bu$ is equipped with the standard grading, then
$L^{-1}N^\bu$ is also equipped with the standard grading
(cf. 4.2). Hence,  
the functors $L$ and $L^{-1}$ induce the functors:
$\DD(L):\DD(\modb(X))\to\DD(\moda'(X))$ and 
$\DD(L^{-1}):\DD(\moda'(X))\to\DD(\modb(X))$
which are quasi-inverse equivalences of triangulated categories.

\subhead{4.3.2}\endsubhead
Recall that $\RR(X)=\RR(X,\text{bottom perversity})$,
$\PP'(X)=\PP(X,\text{top perversity})$.
We will construct a functor 
$\LL^{-1}:\CC(\PP'(X))\to\CC(\RR(X))$. Let 
$\T^\bu\in\CC(\PP'(X))$. We define a bicomplex $\U^{\bu\bu}$
as follows:
$$
\U^{i,j}=\bigoplus_{\dim\D=-j}[\D]^{T^i(\D)}.\tag{$*$}
$$
The differential
$$
\matrix
\U^{i,j}@>{d^{i,j}_I}>>\U^{i,j+1} 
\endmatrix
$$
is the zero map between components $[\D]^{T^i(\D)}$ and 
$[\D']^{T^i(\D')}$
if $\D'\not\sub\D$, and is induced by 
$(-1)^i t^i(\D,\D'):T^i(\D)\to T^i(\D')$
if $\D'$ is a $\op{codim}\ 1$ face of $\D$.\newline
The differential
$$
\matrix
\U^{i,j}@>{d^{i,j}_{II}}>>\U^{i+1,j} \\
\endmatrix 
$$
is the sum of the maps 
$[\D]^{T^i(\D)}\to [\D]^{T^{i+1}(\D)}$ induced by
$d_T^i(\D):T^i(\D)\to T^{i+1}(\D)$.

The diagonal complex of the bicomplex $(*)$ is 
$(\LL^{-1}\T)^\bu$ with the differential 
$d=d_{I}+d_{II}$.

\proclaim{4.3.3. Lemma} The following functorial diagram
$$
\CD
\CC(\PP'(X))@>{\LL^{-1}}>>\CC(\RR(X)) \\
@V{\CC(\Xi_{A'}(X))}VV @V{\CC(\Xi_{B}(X))}VV \\
\CC(\moda'(X))@>{L^{-1}}>>\CC(\modb(X))\\
\endCD
$$
commutes. (Here $\Xi_{A'}(X)=\Xi_{A}(X,\text{top perversity})$.)
\endproclaim

\demo{Proof} The proof follows easily from the fact
that: 
$$
\Hom_k(B(X),N^i_{j})=\bigoplus_{e\in I_{-j}}\II(eN^i)
$$
where $\II(eN^i)$ is the sum of $\dim eN^i$ copies of
the indecomposable injective module 
$\Hom_k(B(X),\F e)$. 
\enddemo

\proclaim{Corollary} The functor 
$\LL^{-1}:\CC(\PP'(X))\to\CC(\RR(X))$
induces a functor 
$\DD(\LL^{-1}):\DD(\PP'(X))\to\DD(\RR(X))$.
\endproclaim

Let us also define a functor $\LL:\CC(\RR(X))\to\CC(\PP'(X))$
as follows: 
$$
\LL=\CC(\Xi_{A'}(X))^{-1}\circ L\circ\CC(\Xi_{B}(X)).
$$
\subhead{4.3.4}\endsubhead
Let $D:\CC(\RR(X))\to\CC(\RR(X))$ be the Verdier duality
functor defined in 2.4. Recall that in 2.4. we also defined
the functors $*$ and $\mp$.

\proclaim{Lemma} The following two functors 
$$
D\simeq\LL^{-1}\circ\CC(\mp)\circ *
$$
are isomorphic.
\endproclaim

\demo{Proof} Let $(\S^\bu,d_S^\bu)\in\CC(\RR(X))$. 
Let $\T^\bu=(\CC(\mp)\circ *)(\S^\bu)^\bu$, and
let $\U^{\bu\bu}$ be a bicomplex defined as in 4.3.2.
Let us also define a bicomplex 
$\V^{i,j}:=D\S^{-i,j}$ (see 2.4 for the definition of 
$D\S^{\bu\bu}$). Then the diagonal complex of $\V^{\bu\bu}$
is precisely $D\S^\bu$. 

Observe that $\V^{i,j}=\U^{i,j}$, and 
$d_{U,I}^{i,j}=(-1)^i d_{V,I}^{i,j}$,
and $(-1)^j d_{U,II}^{i,j}=d_{V,II}^{i,j}$.

Let $q^{i,j}:\V^{i,j}\to\U^{i,j}$ be defined as 
$q^{i,j}=(-1)^{ij}\Id$. Then $q^{\bu,\bu}$ is
an isomorphism of bicomplexes inducing the
isomorphism of the corresponding diagonal complexes.
\enddemo

\proclaim{4.3.5. Corollary} The following two functors
$$
\DD(D)\simeq\DD(\LL^{-1})\circ\DD(\mp)
\circ\DD(*)
$$
are isomorphic.
\endproclaim

\subhead{4.3.6}\endsubhead
Let us simplify the notation for the purposes of this subsection.
We denote $D=\DD(D)$, $L=\DD(\LL)$,
$L^{-1}=\DD(\LL^{-1})$,
$\mp=\DD(\mp)$, and $*=\DD(*)$. We also denote
$\wti{D}=\mp\circ *\circ L^{-1}$.
We call $\wti{D}$ the ``Verdier duality for cosheaves.''

\proclaim{Theorem E} The following functorial diagram
$$
\CD
\DD(\RR(X))@>{D}>>\DD(\RR(X))\\
@V{L}VV @V{L}VV \\
\DD(\PP'(X))@>{\wti{D}}>>\DD(\PP'(X))\\ 
\endCD
$$
commutes. In other words, $L\circ D\simeq\wti{D}\circ L$.
\endproclaim

\demo{Proof} By Corollary 4.3.5:
$$
L\circ D\simeq L\circ L^{-1}\circ\mp\circ *\simeq 
\mp\circ *\simeq\mp\circ *\circ L^{-1}\circ L
\simeq \wti{D}\circ L .
$$
\enddemo

\enddocument
\end